\documentclass[11pt,a4paper]{article}
\pdfoutput=1
\usepackage[T1]{fontenc}
\usepackage[utf8]{inputenc}
\usepackage{lmodern}
\usepackage[english]{babel}
\usepackage[margin=2.4cm]{geometry}
\usepackage{microtype}
\usepackage{amsmath,amssymb,amsthm}
\usepackage{booktabs,array,tabularx}
\usepackage{cancel}
\usepackage{enumitem}
\usepackage[authoryear,round]{natbib}
\usepackage{tikz}
\usetikzlibrary{calc,decorations.pathreplacing,arrows.meta}
\usepackage{xcolor}
\usepackage[hidelinks]{hyperref}
\usepackage{url}

\hypersetup{
  pdftitle={Sharing a unit, measuring a length: a didactic analysis of a teaching progression for fractions from Grade 1 to Grade 8},
  pdfauthor={Vincent Tariel-Adam},
  pdfkeywords={fractions, rational numbers, didactics of mathematics, number line, curriculum design, misconceptions}
}

\newtheorem{principle}{Design principle}
\theoremstyle{definition}
\newtheorem{definition}{Definition}
\newtheorem{hypothesis}{Hypothesis}

\definecolor{fracfill}{RGB}{110,150,200}
\definecolor{fracalt}{RGB}{225,150,90}
\definecolor{fracmix}{RGB}{150,150,150}
\tikzset{
  unitbox/.style={draw,line width=0.9pt},
  partbox/.style={draw,line width=0.4pt},
  shaded/.style={partbox,fill=fracfill},
  shadedalt/.style={partbox,fill=fracalt},
}
\newcommand{\fracbar}[5][shaded]{%
  \pgfmathtruncatemacro{\fbunits}{max(1,ceil(#2/#3))}%
  \pgfmathtruncatemacro{\fbparts}{\fbunits*#3}%
  \foreach \i in {1,...,\fbparts}{%
    \pgfmathsetmacro{\fbx}{(\i-1)*#4/#3}%
    \ifnum\i>#2\relax
      \draw[partbox] (\fbx,0) rectangle ++(#4/#3,#5);%
    \else
      \draw[#1] (\fbx,0) rectangle ++(#4/#3,#5);%
    \fi}%
  \foreach \u in {1,...,\fbunits}{\draw[unitbox] ({(\u-1)*#4},0) rectangle ++(#4,#5);}%
}
\newcommand{\vbar}[5]{%
  \pgfmathsetmacro{\vstep}{2.0/#3}%
  \foreach \i in {1,...,#3}{%
    \pgfmathtruncatemacro{\vk}{\i}%
    \ifnum\vk>#1\relax
      \pgfmathtruncatemacro{\vkk}{\vk-#1}%
      \ifnum\vkk>#2\relax \draw[partbox] (0,{(\i-1)*\vstep}) rectangle ++(0.45,\vstep);%
      \else \draw[partbox,fill=#5] (0,{(\i-1)*\vstep}) rectangle ++(0.45,\vstep); \fi
    \else \draw[partbox,fill=#4] (0,{(\i-1)*\vstep}) rectangle ++(0.45,\vstep); \fi}%
  \draw[unitbox] (0,0) rectangle (0.45,2.0);}
\newcommand{\vbarsub}[4]{%
  \pgfmathsetmacro{\vstep}{2.0/#3}%
  \pgfmathtruncatemacro{\vrest}{#1-#2}%
  \foreach \i in {1,...,#3}{%
    \pgfmathtruncatemacro{\vk}{\i}%
    \ifnum\vk>\vrest\relax
      \ifnum\vk>#1\relax \draw[partbox] (0,{(\i-1)*\vstep}) rectangle ++(0.45,\vstep);%
      \else \draw[partbox,fill=#4!35] (0,{(\i-1)*\vstep}) rectangle ++(0.45,\vstep);
            \draw[partbox] (0,{(\i-1)*\vstep}) -- ++(0.45,\vstep); \draw[partbox] (0.45,{(\i-1)*\vstep}) -- ++(-0.45,\vstep); \fi
    \else \draw[partbox,fill=#4] (0,{(\i-1)*\vstep}) rectangle ++(0.45,\vstep); \fi}%
  \draw[unitbox] (0,0) rectangle (0.45,2.0);}
\newcommand{\numline}[3]{%
  \pgfmathtruncatemacro{\nlticks}{#1*#3}%
  \draw[-{Stealth[length=2mm]}] (-0.2,0) -- ({#3*#2+0.5},0);%
  \foreach \i in {0,...,\nlticks}{\draw ({\i*#2/#1},0.08) -- ({\i*#2/#1},-0.08);}%
  \foreach \k in {0,...,#3}{\draw[thick] ({\k*#2},0.15) -- ({\k*#2},-0.15) node[below=1pt]{\small $\k$};}%
}

\title{Sharing a unit, measuring a length:\\ a didactic analysis of a teaching progression\\ for fractions from Grade~1 to Grade~8}
\author{Vincent Tariel-Adam\\[2pt]
  \small Independent researcher; CommeUnJeu (\url{https://commeunjeu.org}), France\\
  \small \texttt{vincent.tariel@alumni.polytechnique.org}}
\date{October 2026}

\begin{document}
\maketitle

\begin{abstract}
\noindent
Fractions require learners to revise properties associated with whole numbers. This paper presents a design rationale and a critical analysis of an eight-level source progression for fractions (Grades 1--8, nominal ages 6--14), from which the French-curriculum and International Baccalaureate chapters published at \url{https://commeunjeu.org} are assembled. Its central hypothesis is a bridge from sharing to numerical measure: a share of one unit is represented by a bar, read as a length relative to a specified unit, and located on the number line. The reference unit is kept stable in the initial tasks to focus attention on fraction magnitude, equivalence and operations; changing the unit is deliberately postponed, without treating this as a mathematical requirement or an established reduction in cognitive load. Three principles organise the meanings: sharing one unit, coordinating bars and number lines, and establishing fractions as numbers before explicitly treating them as scalar operators. Two supporting principles concern valid local justifications before formal generalisation and chains of equalities, introduced in Grades 2--3, with written reasons in lessons from Grade 4. The analysis relates these choices to the research literature and selected curricula, examines task strategies and feedback, and distinguishes realised features from discrepancies and proposed revisions. The source documents and their curriculum projections do not establish the order actually experienced by pupils. No empirical evaluation is reported. Primary hypotheses address transfer between representations, coordination of units and the numerical status of fractions before operators; the audit also identifies limits to the reproducibility of its findings.

\end{abstract}


% ======================================================================
\section{Introduction}\label{sec:intro}
% ======================================================================

Fractions occupy a singular place in school mathematics. Extending whole-number arithmetic to rational numbers requires learners to reconsider notation, order and operations: distinct fractional numerals such as $\frac{1}{2}$, $\frac{2}{4}$ and $\frac{3}{6}$ denote the same number; no rational number has an immediate successor in the ordered set of rational numbers; multiplying a positive quantity by a number strictly between $0$ and $1$ makes it smaller; and, for a fixed unit and a fixed positive numerator, increasing a positive denominator decreases the fraction's magnitude. Psychological research describes the systematic transfer of whole-number properties to fractions as a \emph{whole-number bias} \citep{ni2005}. Longitudinal studies have found that knowledge of fractions and of division at age ten predicts mathematics and algebra achievement in high school, beyond whole-number arithmetic, general cognitive measures and family background \citep{siegler2012}, that fraction knowledge predicts gains in achievement from one year to the next \citep{bailey2012}, and that the relation holds in samples from three continents \citep{torbeyns2015}; fraction magnitude knowledge is also related to readiness for algebra \citep{booth2012}. These findings are correlational, but they have made fractions a recognised priority \citep{nmap2008}. Fractions are, in this sense, a turning point of the curriculum rather than one chapter among others.

Mathematics-education research has long described the rational number as a family of related constructs---part--whole, measure, quotient, ratio and operator \citep{kieren1976,behr1983,lamon2007}---each mathematically legitimate, and each making different demands on the learner. A curriculum author therefore faces a design problem that is not to choose the correct meaning, since all of them are correct, but to choose a \emph{first} meaning that can carry the others without being contradicted by them, a \emph{representation} that can accompany that meaning over several years, and an \emph{order} in which the other meanings are added.

\begin{figure}[htbp]
\centering
\begin{tikzpicture}[
  stage/.style={draw=fracfill!65!black,line width=0.6pt,rounded corners=2pt,
    text width=3.15cm,minimum height=3.45cm,inner sep=6pt,align=center,
    font=\small},
  transition/.style={-{Stealth[length=2.3mm]},line width=0.8pt,
    draw=fracfill!65!black}
]
\node[stage,fill=fracfill!10] (sharing) at (0,0) {
  \begin{minipage}[t][2.9cm][t]{3.15cm}\centering
  \textbf{1. Sharing\\one unit}\\[8pt]
  Constructing unit\\
  fractions and\\
  iterating them
  \end{minipage}
};
\node[stage,fill=fracfill!14] (measure) at (4.05,0) {
  \begin{minipage}[t][2.9cm][t]{3.15cm}\centering
  \textbf{2. Measuring\\a length}\\[8pt]
  Using the bar as\\
  an intermediate\\
  representation
  \end{minipage}
};
\node[stage,fill=fracfill!18] (number) at (8.1,0) {
  \begin{minipage}[t][2.9cm][t]{3.15cm}\centering
  \textbf{3. Establishing\\number status}\\[8pt]
  Number line, order,\\
  equivalence and\\
  arithmetic
  \end{minipage}
};
\node[stage,draw=fracalt!65!black,fill=fracalt!12] (extensions) at (12.15,0) {
  \begin{minipage}[t][2.9cm][t]{3.15cm}\centering
  \textbf{4. Extending uses\\and definitions}\\[8pt]
  Operator, quotient,\\
  signed rational\\
  numbers and\\
  algebraic fractions
  \end{minipage}
};
\draw[transition] ([xshift=2pt]sharing.east) -- ([xshift=-2pt]measure.west);
\draw[transition] ([xshift=2pt]measure.east) -- ([xshift=-2pt]number.west);
\draw[transition] ([xshift=2pt]number.east) -- ([xshift=-2pt]extensions.west);
\end{tikzpicture}
\caption{Four stages of the proposed progression. Arrows indicate the intended order of explicit teaching. The bar provides an intermediate representation between sharing a unit and establishing fractions as numbers through length measurement. The initial stages use a stable reference unit; early sharing tasks may already involve uses developed explicitly later.}
\label{fig:overview}
\end{figure}
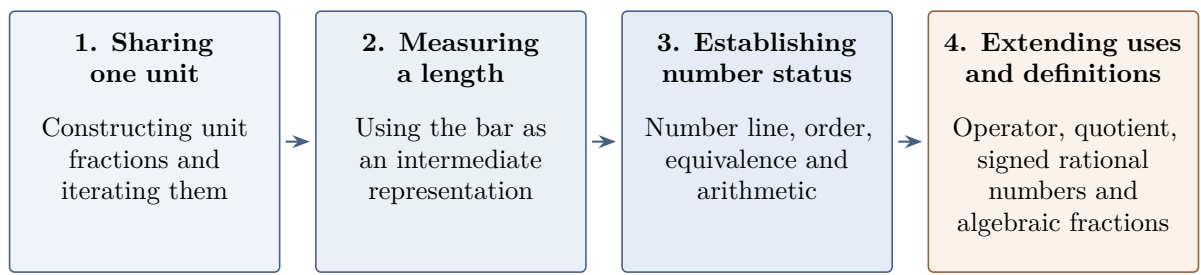

This paper documents and analyses one answer to that design problem: an eight-year progression for fractions, from Grade~1 to Grade~8 (pupils aged about 6 to 14), written as bilingual (English and French) lessons and exercise sets. Its thesis is that learning fractions means crossing a bridge between two worlds---the world of \emph{sharing}, in which a thing is cut fairly and parts are counted, and the world of \emph{measure}, in which a length is compared with a unit and expressed by a number---and that the design of the progression is organised around that crossing. Figure~\ref{fig:overview} summarises the proposed order of explicit teaching in four stages. The design rests on five principles:
\begin{enumerate}[label=(\roman*)]
  \item a single foundational definition in the world of sharing: the fraction as a fair share of one continuous unit, whose reading is extended until the quotient takes over as the operative definition (Section~\ref{sec:p1});
  \item the bar as the bridge to the world of measure, and the number line as the arrival: the fraction is represented successively by a measured length and a point, with the aim of establishing its numerical status (Section~\ref{sec:p2});
  \item the numerical status of fractions is established before their explicit treatment as scalar operators; applications to quantities and part--part ratios are developed in a companion strand (Section~\ref{sec:p3});
  \item technique before formalism: tools such as cross-multiplication should receive a valid local justification and age-appropriate validity conditions when introduced, with formal generalisation deferred; the audit identifies where the chapters fall short (Section~\ref{sec:p4});
  \item a common written form for multi-step calculations, the chain of equalities, introduced in the primary grades and annotated in lessons from Grade~4, with the intended extension to algebra examined separately (Section~\ref{sec:p5}).
\end{enumerate}

Principles~1 to~3 form the progression of meanings and are the paper's main argument; Principles~4 and~5 are supporting instructional decisions, kept because they shape how the chapters are written. Section~\ref{sec:progression} describes the progression grade by grade, Section~\ref{sec:misconceptions} analyses the misconceptions it anticipates and the didactic variables of its central task, Section~\ref{sec:implementation} describes how it is implemented, and Section~\ref{sec:discussion} confronts the design with the chapters, discusses its limitations and states the hypotheses that an empirical evaluation should test.

\paragraph{Genre and scope.} The paper is a design rationale and a curriculum analysis. Like a hypothetical learning trajectory \citep{simon1995} or a research-based learning trajectory \citep{clements2004}, the progression embodies learning goals, an order of activities and hypotheses about how learning proceeds; the paper makes these explicit and confronts them with the literature and with the chapters that implement them. Its method borrows from the \emph{a priori analysis} of didactic engineering \citep{artigue1990,artigue1992} at the level of global, or macro-didactic, choices---the definition, the representation, the order of meanings---and Section~\ref{sec:taskvariables} gives one example at the level of a task and its didactic variables. It does not contain the \emph{a posteriori} part of an engineering study: no classroom experiment or controlled comparison is reported, and statements about learners are either taken from the cited literature or presented as design hypotheses. None of the five principles is new taken alone; several are close to recommendations already present in research syntheses and curricula \citep{siegler2010,ccss2010}. What the paper contributes is their articulation across eight source levels, the explicit statement of what each choice postpones, and a confrontation of the stated design with the materials that implement it, which brings out where the two diverge (Section~\ref{sec:discussion}). The principles have a dual role: they express the intended design and provide criteria for revising its current realisation. Descriptions of existing blocks, interpretations of their didactic role and proposed changes are therefore distinguished throughout; the article does not claim that every principle is already fully implemented.

\paragraph{Corpus and method.} The corpus is the sixteen source documents of the fraction progression (a lesson and an exercise set for each of Grades~1 to~8) and the five lessons of the companion ratio strand (Grades~2 to~6), in the version of 9~October~2026. The unit of analysis is the \emph{block}: the chapters are written in a mark-up in which every definition, proposition, method, example, exercise and answer is a delimited block with a type, so that statements such as ``the Grade~3 lesson states the general procedure'' or ``three of the thirteen simplification statements with a letter in the denominator give the condition'' are made by locating and counting blocks of a given type, not by impression. Three kinds of statement are made about the corpus: (i) \emph{presence}, that a block with a given content exists at a given grade; (ii) \emph{order}, that one block precedes another within a chapter or across grades; (iii) \emph{counts}, which are exhaustive over the mark-up (for instance, the number of multi-line calculations and, among them, the number carrying a written reason). Ambiguous cases are reported as such (for example, a condition stated in the sentence before a display rather than beside a line). 

For the counts of Principle~5, a \emph{chain} is a multi-line calculation, which the mark-up always writes in one environment (\texttt{align} or \texttt{aligned}), and a chain carries a \emph{reason} when one of its lines bears a parenthesised annotation such as ``(divide both by $2$)''; words that belong to the expression itself (``of'', ``number of girls'') do not count. Every lesson and exercise is written twice in the source, in English and in French, and the counts are made on the English half only. They were made by a script over the mark-up and checked by reading. This criterion records the presence of reasons, not the justification of every step that Principle~5 describes, and gives an upper bound on the chains that fully realise it; under the stricter criterion of a reason on every line except the last, which usually records the result, the lessons of Grades~4 to~8 annotate $30$ of their $54$ chains, and the exercise answers $48$ of $477$. 

Every statement was checked against the source files after the principles had been formulated, and the discrepancies found are reported in Section~\ref{sec:discussion}; Table~\ref{tab:trace} gives, for each principle, the blocks that carry it. The counts concern the source documents, not repeated appearances of a block in assembled curriculum chapters. Single-line equality chains and reasons given only in surrounding prose fall outside this annotation count; their exclusion does not establish that a calculation is unjustified. Lessons and exercise answers are reported separately. Appendix~\ref{app:audit} summarises these conventions and distinguishes them from requirements for an independently inspectable reanalysis. The source documents and the counting script are held by the author and can be consulted on request; the assembled curriculum chapters are public (Appendix~\ref{app:resources}). No public archive of the dated source corpus, script and block-level coding accompanies this manuscript, so the reported counts are not yet independently reproducible from the public chapters alone.

\paragraph{Grade convention.} The progression is organised in eight levels, labelled Grade~1 to Grade~8 and nominally corresponding to ages 6 to 14 (in the French system, CP to \emph{quatri\`eme}). It is a \emph{source} progression rather than a national curriculum: the chapters of a given curriculum are assembled from it topic by topic, so that a level of the progression is not bound to one school class. In the French-curriculum chapters, for instance, fractions start in CE1 (age 7), and the content of Grade~6 is distributed over three classes (Section~\ref{sec:implementation}). Unless stated otherwise, ``Grade~$N$'' refers to the levels of the progression.

% ======================================================================
\section{Theoretical background}\label{sec:background}
% ======================================================================

\subsection{One object, several constructs}

\citet{kieren1976} distinguished several interpretations of rational numbers, among them measures and points on a number line, which he later organised into a small number of \emph{sub-constructs} with partitioning and equivalence as their mechanisms \citep{kieren1980}. \citet{behr1983} proposed a model in which the part--whole construct, supported by partitioning, is fundamental to the ratio, operator, quotient and measure constructs. \citet{freudenthal1983}, in his didactical phenomenology, analysed the fraction as what fractures, what compares and what operates. \citet{lamon2007} treats fractions, ratios and proportions as one multiplicative domain, among the most difficult to teach, and proposes a framework for research on it. A test of the model of \citet{behr1983} with Grade~5 and~6 pupils \citep{charalambous2007} found the part--whole construct to be necessary but not sufficient for the others, and only weakly related to the quotient and measure constructs; the authors note that this hierarchy may reflect a curriculum in which part--whole is taught first.

A recurrent criticism concerns the \emph{form} that the part--whole meaning takes in teaching. When $\frac{a}{b}$ is read as ``$a$ parts out of $b$'', children count shaded and unshaded pieces, and the reading breaks down for improper fractions such as $\frac{5}{3}$, where there are not five parts ``out of'' three; the limits of the part-of-a-whole model were already a central finding of studies of secondary pupils' errors \citep{kerslake1986}. Studies of children's fraction schemes \citep{steffe2010} distinguish several readings. In a \emph{part--whole} scheme, $\frac{a}{b}$ is $a$ of the $b$ parts of a partitioned whole. In the \emph{partitive} schemes, $\frac{1}{b}$ is a part which, iterated $b$ times, reproduces the whole, and $\frac{a}{b}$ is $a$ iterations of $\frac{1}{b}$, still contained in the whole. In the \emph{iterative fraction scheme}, the unit fraction is freed from its containment in the whole, so that $a$ may exceed $b$ and $\frac{5}{3}$ becomes a number. The passage to the iterative scheme requires a reorganisation of the child's units \citep{tzur1999,hackenberg2007}; reading $\frac{m}{n}$ as ``$m$ one-$n$-ths'' rather than ``$m$ out of $n$'' \citep{thompson2003} expresses the same idea. The distinction matters for Section~\ref{sec:p1}: ``parts of a unit'' can be taught in any of these forms.

\subsection{Fractions as numbers with magnitude}

The \emph{integrated theory of numerical development} \citep{siegler2011,siegler2013} describes numerical development as coming to understand that all real numbers, and not only whole numbers, have magnitudes that can be ordered and located on number lines; on this view, what fractions share with whole numbers is not their notation but the fact that they have magnitudes. The number line is, in that theory, a way of assessing magnitude knowledge as well as a means of teaching it; the theory does not by itself prescribe a representation. The practice guide of \citet{siegler2010} accordingly recommends helping students to recognise that fractions are numbers that expand the number system beyond whole numbers, and to use number lines as a central representational tool; it also recommends building on students' informal understanding of sharing, while rating the evidence for that recommendation as minimal and the evidence for number lines as moderate. Intervention studies that organise fraction teaching around magnitude and the number line \citep{fuchs2013,hamdan2017} have reported gains in fraction magnitude understanding over control conditions, and \citet{gunderson2019} attribute part of the advantage of the number line to its unidimensionality. Earlier design work had already explored linear and measurement entries into rational numbers \citep{moss1999}.

\subsection{Conceptual fields}

\citet{vergnaud1983,vergnaud1990,vergnaud2009} defines a concept by a triplet: the set of situations that give it meaning, the operational invariants that the learner uses to deal with those situations, and the linguistic and symbolic representations that express them. Fractions belong to the conceptual field of \emph{multiplicative structures}, together with multiplication, division, ratio, proportion and linear functions; the prototypical situation of that field is the \emph{isomorphism of measures}, a proportional relation between two measure spaces. In these terms, the design problem of the introduction becomes: in which order should the situations of the multiplicative field be attached to the fraction notation, and which invariants should be secured before new situations are added?

\subsection{The French didactic tradition}

The theory of didactical situations \citep{brousseau1997} analyses the conditions under which a piece of knowledge appears as the solution to a problem that previous knowledge cannot solve. Its best-known application to our topic is the long engineering on rationals and decimals \citep{brousseau1987,brousseau2004,brousseau2007,brousseau2008,brousseau2009}, in which pupils first meet rational numbers as measures obtained by \emph{commensuration}---designating the thickness of a sheet of paper by a pair of numbers, a number of sheets and the thickness of the pile---and later as decimals and as linear maps. In the same research programme, \citet{ratsimba1981} compared two methods of rational measurement, \emph{commensuration} and the \emph{fractionnement de l'unit\'e} (partitioning the unit); the progression analysed here uses the second, and the contrast is discussed in Sections~\ref{sec:p1} and~\ref{sec:discussion}. The theory of \emph{didactic transposition} \citep{chevallard1985,chevallard1991} describes the transformations that a piece of scholarly knowledge undergoes when it becomes knowledge to be taught; together with the praxeological description of the anthropological theory of the didactic \citep{chevallard1999}---a type of task, a technique, a technology that justifies the technique, a theory that justifies the technology---we use it in Section~\ref{sec:p4} to describe what is and is not transposed when cross-multiplication is taught.

\subsection{Institutional context}

The progression was written to be compatible with the French national programmes, which it covers. Until 2025 those programmes introduced fractions in CM1 \citep{men2020programmes}; the programmes in force since 2025 introduce them in CE1 as parts of a whole---halves and quarters first, then unit fractions, then other fractions not greater than one---and in CE2 place them on a graduated unit strip, an approach which, in the programme's words, helps ``to free [pupils] from the `whole' and give fractions the status of numbers'' \citep{men2024cycle2}. In CM1, fractions greater than one are placed on a graduated half-line, the unit fraction being treated as a new unit of measure; in \emph{sixi\`eme} the fraction acquires the meaning of a quotient \citep{men2025cycle3}, which national guidance verbalises as ``the quotient of $a$ by $b$ is the number which, multiplied by $b$, gives $a$'' \citep{men2025exemples6e}. The French programmes thus follow the same order of meanings as the progression analysed here: part of a unit, then measure and point on a line, then quotient; Section~\ref{sec:p3} examines the one meaning on which they differ.
Its order is also close to that of the Common Core State Standards \citep{ccss2010}, where unit fractions, fractions on the number line and equivalence ``as the same point on a number line'' appear in Grade~3; the general equivalence rule, addition with like denominators, the product of a fraction by a whole number and decimal fractions in Grade~4; and unlike denominators, the product of two fractions (with an area model) and the quotient interpretation in Grade~5. The main differences with both are discussed in Section~\ref{sec:discussion}.

% ======================================================================
\section{Principle 1: a fair share of one unit}\label{sec:p1}
% ======================================================================

\begin{principle}[A fair share of one unit]\label{pr:unit}
A fraction is introduced as the result of sharing one unit fairly: the unit is cut into $b$ equal parts and $a$ of them are taken. The unit is one continuous object, never a collection. This definition belongs to the world of sharing; Principle~\ref{pr:linear} is intended to establish its numerical interpretation through measure, and from Grade~6 the quotient characterisation ($\frac{a}{b}$ is the number which, multiplied by $b$, gives $a$) becomes the operative definition.
\end{principle}

\paragraph{Terminology.} The distinction between \emph{unit} and \emph{whole} is a convention of this design, not a general restriction on mathematical usage. In the world of sharing, \emph{unit} names the continuous reference object whose equal parts are considered; \emph{whole} is used for a collection or quantity such as twelve beads. In the world of measure, a reference length $U$ is assigned measure~$1$, and a length $L$ has numerical measure $r$ when $L=rU$. The number~$1$, the reference length and its physical or graphic representation are distinct. A paper strip has a physical length even when no value in centimetres is specified. In early tasks the reference is kept stable, and ``a length of $\frac{2}{3}$'' abbreviates a length measuring $\frac{2}{3}$ in that reference unit. Changing the reference is mathematically legitimate; the reasons for postponing it are didactic (Section~\ref{sec:crossing}).

\subsection{Three candidate entry points}

\begin{figure}[ht]
\centering
\begin{tikzpicture}
\begin{scope}
\fracbar{2}{3}{3.9}{0.6}
\draw[decorate,decoration={brace,amplitude=4pt,mirror}] (0,-0.1) -- (2.6,-0.1) node[midway,below=5pt]{2 of the 3 equal parts};
\node[above] at (1.95,0.65) {\small one unit, cut fairly in 3};
\node at (1.95,-1.6) {(a) a share of one unit};
\end{scope}
\begin{scope}[xshift=5.3cm]
\foreach \i in {0,...,5}{
  \pgfmathsetmacro{\cx}{0.55+\i*0.62}
  \ifnum\i<4 \draw[fill=fracalt] (\cx,0.3) circle (0.24); \else \draw (\cx,0.3) circle (0.24); \fi}
\node[below] at (2.1,-0.15) {$\tfrac{4}{6}=\tfrac{2}{3}$ of the apples};
\node[above] at (2.1,0.65) {\small 4 red apples out of 6};
\node at (2.1,-1.6) {(b) a part of a collection};
\end{scope}
\begin{scope}[xshift=10.7cm]
\begin{scope}[yshift=0.35cm]\fracbar{1}{3}{1.8}{0.35}\end{scope}
\begin{scope}[xshift=2.1cm,yshift=0.35cm]\fracbar{1}{3}{1.8}{0.35}\end{scope}
\node[above] at (1.95,0.75) {\small 2 units shared among 3};
\begin{scope}[yshift=-0.45cm,xshift=0.9cm]\fracbar{2}{3}{1.8}{0.35}\end{scope}
\node[right] at (2.8,-0.28) {\small one share};
\node[below] at (1.8,-0.55) {$2\div 3=\tfrac{2}{3}$};
\node at (1.95,-1.6) {(c) a quotient};
\end{scope}
\end{tikzpicture}
\caption{Three entry points for $\frac{2}{3}$. (a) One unit is cut fairly into three parts and two of them are taken. (b) Four of the six apples are red: the fraction describes a relation between a sub-collection and a collection. (c) Two units are shared equally among three people: each share consists of one third of each unit, that is, two thirds of a unit.}
\label{fig:entry}
\end{figure}
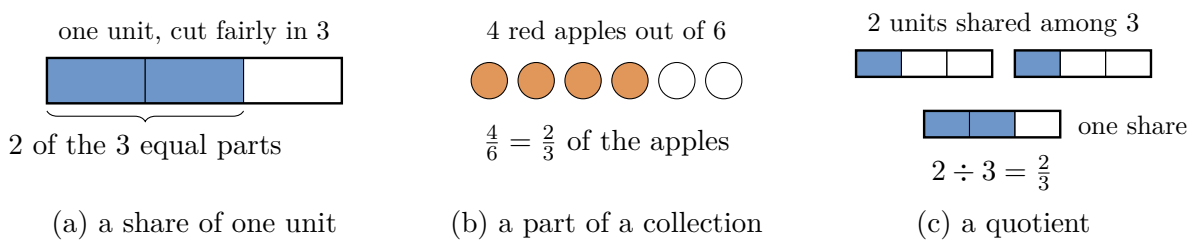

Figure~\ref{fig:entry} shows three entry points that a first lesson on $\frac{2}{3}$ can realistically take, and Table~\ref{tab:entry} compares them from the designer's point of view.

\begin{table}[ht]
\centering\small
\begin{tabularx}{\textwidth}{>{\raggedright\arraybackslash}p{2.5cm}>{\raggedright\arraybackslash}X>{\raggedright\arraybackslash}X>{\raggedright\arraybackslash}X}
\toprule
 & (a) Share of one unit & (b) Part of a collection & (c) Quotient \\
\midrule
Referent & one fixed continuous unit & a collection, which changes from one situation to the next & $a$ units shared among $b$ sharers \\
What $\frac{a}{b}$ names & a quantity: $a$ parts, each one of $b$ equal parts of the unit & a relation between a sub-collection and the collection & the size of one share when $a$ units are shared among $b$ \\
Actions required & cut one unit fairly; count parts & count two collections; relate them & cut each unit; coordinate $a$ units and $b$ sharers; recombine one share \\
Passage to measure & immediate on a bar: every fair cut of a bar is a length & only through an added representation (equal objects drawn as equal lengths) & after recombination \\
Natural continuation & number line; fractions greater than one & proportion; percentages; the operator & division; the field of rational numbers \\
Main risk as a first meaning & ``$a$ out of $b$'' counting, which does not extend beyond one & confusing the fraction with a cardinal; part--part versus part--whole & coordinating $a$ units and $b$ sharers on a first encounter \\
\bottomrule
\end{tabularx}
\caption{A priori comparison, from the designer's point of view, of the three entry points of Figure~\ref{fig:entry}.}
\label{tab:entry}
\end{table}

The comparison is not a ranking of mathematical value; it states, as design hypotheses, what each entry asks a young learner to coordinate on a first encounter, and whether it opens onto measure. Entries (a) and (c) are closer than they look: the unit fraction $\frac{1}{b}$ is itself the share obtained when one unit is shared among $b$, that is, the quotient $1\div b$, and under (a) the fraction $\frac{a}{b}$ is $a$ such shares---the first situation of the progression is indeed the fair sharing of one pizza. The difference lies in the number of units: (c) requires coordinating $a$ units and $b$ sharers, and the share it produces is made of parts taken from different units. Entry (b) differs in its referent: a collection has no length of its own, so the passage to measure, which (a) offers at once on a bar, requires an added representation. Research on equal sharing shows that sharing situations with several units can supply the meaning of fractions and that children solve them with informal partitions before they meet the notation \citep{streefland1991,empson2011}; naming one share in two ways, for instance $\frac{1}{2}+\frac{1}{6}$ and $\frac{2}{3}$, is then a natural source of equivalence, and pupils' informal knowledge of parts can serve as a basis for the symbols \citep{mack1990}. The choice of (a) is therefore a choice of order---one unit before several, and sharing before measure---not a claim that the other entries are inaccessible.

\subsection{How the progression enters}

The first lesson of Grade~1 does not begin with the notation but with the condition on which sharing rests. A pizza is cut in two different ways, and the pupil decides which cut is fair. Halves and quarters are then named and written $\frac{1}{2}$ and $\frac{1}{4}$, on discs, rectangles and strips, by folding a strip of paper, and on collections, where half of twelve strawberries is found by sharing them equally into two groups; relations such as ``two quarters make a half'' are stated in words. The general notation $\frac{a}{b}$, with its numerator and denominator, comes afterwards, with a single situation---one cake cut into two equal parts, of which one is taken---and a definition read on a bar:

\begin{definition}[Fraction, Grade~1, in substance]\label{def:fraction}
A fraction has two numbers separated by a fraction bar. The \emph{denominator}, below the bar, tells how many equal parts the unit is divided into; the \emph{numerator}, above the bar, tells how many of those parts are taken.
\end{definition}

\begin{figure}[ht]
\centering
\begin{tikzpicture}[>={Stealth[length=2mm]}]
  % the lesson's figure: a vertical unit bar cut into three, two parts shaded from the bottom,
  % an equals sign, the numeral, and the two numbers coloured like what they count
  \begin{scope}[shift={(-1.7,-1.0)}]
    \foreach \i in {1,2}{\draw[partbox,draw=fracalt,fill=fracfill] (0,{(\i-1)*2.0/3}) rectangle ++(0.45,2.0/3);}
    \draw[partbox,draw=fracalt] (0,{2*2.0/3}) rectangle ++(0.45,2.0/3);
    \draw[unitbox] (0,0) rectangle (0.45,2.0);
  \end{scope}
  \node at (-0.8,0) {$=$};
  \node[font=\Large,text=fracfill] (num) at (0,0.42) {$2$};
  \node[font=\Large,text=fracalt] (den) at (0,-0.45) {$3$};
  \draw[thick] (-0.22,0) -- (0.22,0);
  \draw[thick,<-] ([xshift=2pt]num.east) -- ++(0.9,0) node[right,text width=6.4cm,align=left]{The \textbf{numerator} tells how many parts we have (shaded).};
  \draw[thick,<-] ([xshift=2pt]den.east) -- ++(0.9,0) node[right,text width=6.4cm,align=left]{The \textbf{denominator} tells how many equal parts the unit is divided into.};
\end{tikzpicture}
\caption{The Grade~1 lesson's figure for the definition. The unit bar, cut into three equal parts with two of them shaded, stands beside the numeral, joined by an equals sign; each number of the fraction is coloured like what it counts, the numerator like the shaded parts and the denominator like the cuts of the unit. The bar is drawn vertically, as in the lessons' addition figures (Figure~\ref{fig:operations}).}
\label{fig:def}
\end{figure}
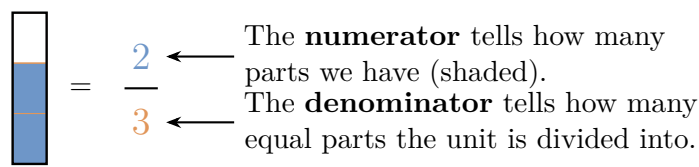

The lesson does not leave the definition to the words: Figure~\ref{fig:def} reproduces its figure, in which the same object is written twice---as a bar and as a numeral---and the equals sign between them is the first statement that the drawing and the notation name one thing. The definition takes the unit, not a collection, as its referent, and it requires the parts to be equal: fairness is the condition under which a share can later be read as a measure, which is why the chapter opens on a fair and an unfair cut. Partitioning into equal parts is itself a competence that develops in stages in young children \citep{pothier1983}, which is one reason the first grade restricts itself to halves and quarters before the general notation. Since the denominator counts the parts into which the unit is cut, it is a whole number at least equal to~$1$; the Grade~1 chapter leaves this implicit, and the condition is stated in words in Grade~4 (Section~\ref{sec:p4}). As worded, the definition is the ``$a$ out of $b$'' counting reading discussed in Section~\ref{sec:background}, which is unproblematic as long as $a\leq b$. The progression then extends that reading in four steps, each of which keeps the definition and adds one property:
\begin{enumerate}[label=(\arabic*)]
  \item Grade~2 defines the \emph{unit fraction} as one single part and states that every fraction is made of unit fractions: $\frac{3}{8}=\frac{1}{8}+\frac{1}{8}+\frac{1}{8}$, read ``three eighths is three times one eighth''. This is the partitive reading of \citet{steffe2010}---$\frac{a}{b}$ as $a$ iterations of $\frac{1}{b}$, still contained in the unit---and it coincides with the definition of the Common Core standards, where $\frac{a}{b}$ is the quantity formed by $a$ parts of size $\frac{1}{b}$ \citep{ccss2010}, and it is compatible with the definition of a fraction as the $m$-th point in the sequence of $n$-ths on the number line \citep{wu2011}, which Grade~3 will supply. Writing fractions as sums of unit fractions has a long history: Egyptian scribes expressed fractions as sums of distinct unit fractions \citep{imhausen2016}; the Grade~2 decomposition repeats one unit fraction instead, which is what makes it a bridge to the product $a\times\frac{1}{b}$.
  \item Grade~4 lets the counting continue past the unit, every successive unit being partitioned in the same way: $\frac{3}{2}$ is three halves, one unit and half of a second one, and $\frac{5}{3}=1+\frac{2}{3}$ is placed on the number line beyond~$1$. Strictly, this extends the definition from parts of one unit to parts of the same size taken from as many units as needed, the point at which the Grade~1 wording no longer suffices; it is also where the iterative fraction scheme is called for. It is a step that the world of measure makes natural---a ruler does not stop at one unit---and the world of sharing does not.
  \item Grade~6 writes the iteration as a product, $\frac{n}{b}=n\times\frac{1}{b}$, and relates proper fractions, improper fractions and mixed numbers.
  \item Grade~7 extends the numbers to signed fractions. Here the sharing reading cannot be extended (there is no negative number of parts), and the progression relies instead on the quotient meaning of Grade~6: the fraction bar is read as a division of integers, so that the sign rules of division give $\frac{-a}{b}=\frac{a}{-b}=-\frac{a}{b}$; the opposite of $\frac{a}{b}$ is $-\frac{a}{b}$, since their sum is $0$; and a fraction is rewritten with a positive denominator before it is compared or added.
\end{enumerate}
The partitive reading is thus prepared two years before fractions greater than one call for the iterative one, the transition at which the counting reading is known to fail \citep{tzur1999,thompson2003} and which the French cycle~3 programme also identifies \citep{men2025cycle3}. Whether a textbook statement suffices to bring about the reorganisation of units that this transition requires \citep{tzur1999,hackenberg2007} is an empirical question (Hypothesis~\ref{h:improper}).

Steps (3) and (4) mark a hand-over. In Grade~6 the chapter infers from the quotient, $a\div b=\frac{a}{b}$, that $\frac{a}{b}$ is the number which, multiplied by $b$, gives $a$---the characterisation that the French programmes adopt in \emph{sixi\`eme} \citep{men2025exemples6e}---and Grades~7 and~8 treat the fraction bar as a division sign. The sharing definition therefore remains the definition of non-negative fractions as quantities, while the quotient becomes the definition used for computation and for signed and algebraic fractions. Principle~\ref{pr:unit} is best read in that form: one definition carries the meaning until a second one, introduced as a property of numbers already known, takes over the computations. How the hand-over is justified is a question for Principle~\ref{pr:technique} (Section~\ref{sec:p4}).

Collections are not absent from the first grades. Half or a quarter of a collection (Grade~1) and a third of twelve beads (Grade~2) appear as applications of fair sharing to a whole: the fraction names a share of a quantity. What the first grades do not do is define the fraction on a collection, or use a collection as the referent of the notation; Section~\ref{sec:p3} explains where ``$\frac{a}{b}$ of a quantity'' is developed.

\subsection{What the choice gains and what it postpones}

The choice has four didactic motivations. (1) The analysis of Table~\ref{tab:entry} suggests that it may reduce the number of relations to coordinate on a first encounter; no comparative cognitive-load measurement is available. (2) On a bar, the share is a length from the first day, which is what Principle~\ref{pr:linear} needs. (3) It allows the other meanings to be added later as enrichments rather than as corrections. (4) It lends itself to the Concrete--Pictorial--Abstract sequence, an adaptation \citep{leong2015} of Bruner's enactive, iconic and symbolic modes of representation \citep{bruner1964,bruner1966}: strips that are folded and cut, bars that are drawn and shaded, then the notation $\frac{a}{b}$.

The choice also postpones three things, for which the progression pays later. The quotient meaning, which some research programmes put first through equal-sharing problems with several objects \citep{streefland1991,empson2011}, waits until Grade~6, where two cakes are shared among three people and $a\div b=\frac{a}{b}$ is stated. The written operator, ``$\frac{a}{b}$ of $N$'' as the multiplication $\frac{a}{b}\times N$, is kept out of the fraction strand until the same grade (Section~\ref{sec:p3}). And the commensuration entry of \citet{brousseau2004}, in which a fraction first designates a relation between two measurable quantities, is not used: of the two methods of rational measurement compared by \citet{ratsimba1981}, the progression uses only the partitioning of a given unit. Section~\ref{sec:discussion} returns to the cost of these choices.

% ======================================================================
\section{Principle 2: the bridge from sharing to measuring}\label{sec:p2}
% ======================================================================

\begin{principle}[The bar as a bridge, the number line as arrival]\label{pr:linear}
The share is represented by a bar cut into equal parts: an object that is at once a shared unit and a length. The bar is intended to support a transition from sharing to measure. In Grade~3, rulers graduated in fractions of a specified unit and points on the number line explicitly represent the fraction as a number. This describes the teaching intention, not evidence that the transition has occurred for each pupil. Other representations are used for specific purposes and related back to the bar.
\end{principle}

\subsection{Two worlds and one object}

Fair sharing and measuring are two different activities. In the first, a thing is cut and parts are counted; the result is ``two of the three parts of the cake''. In the second, a length is compared with the unit, the length taken as~$1$, and the result is a number, ``$\frac{2}{3}$'', which can be placed, compared and added like other numbers. \citet{freudenthal1983} distinguishes the fraction as \emph{fracturer}---what breaks a whole into parts---from the fraction as \emph{comparer}, which relates two magnitudes; the bridge described here leads from the first to one case of the second, comparison with the unit length~$1$. In the sub-construct model, partitioning is the mechanism and measure the construct it leads to \citep{kieren1980,behr1983}; and the French programmes describe the same passage, in the same year as this progression, when the CE2 pupil graduates a unit strip, an approach meant ``to free [pupils] from the `whole' and give fractions the status of numbers'' \citep{men2024cycle2}.

\begin{figure}[ht]
\centering
\begin{tikzpicture}
\begin{scope}
\fracbar{2}{3}{6}{0.5}
\draw[decorate,decoration={brace,amplitude=4pt}] (0,0.6) -- (4,0.6) node[midway,above=5pt]{\small sharing: 2 of the 3 equal parts of the unit};
\draw[decorate,decoration={brace,amplitude=4pt,mirror}] (0,-0.1) -- (4,-0.1) node[midway,below=5pt]{\small measure: a length of $\tfrac{2}{3}$};
\end{scope}
\draw[-{Stealth[length=2.5mm]},thick] (6.6,0.25) -- (7.6,0.25) node[midway,above=2pt]{\scriptsize Grade 3};
\begin{scope}[xshift=8.2cm,yshift=0.25cm]
\numline{3}{4.5}{1}
\fill[fracfill!60!black] (3,0) circle (2.2pt);
\node[above=3pt] at (3,0) {$\tfrac{2}{3}$};
\node[below=12pt] at (2.25,0) {\small number: a point};
\end{scope}
\end{tikzpicture}
\caption{The bridge. One bar, two readings: a shared unit (2 of the 3 equal parts) and a length ($\frac{2}{3}$, the bar itself having length $1$). Placed from $0$, the measured segment ends at the point representing the fraction on the number line.}
\label{fig:bridge}
\end{figure}
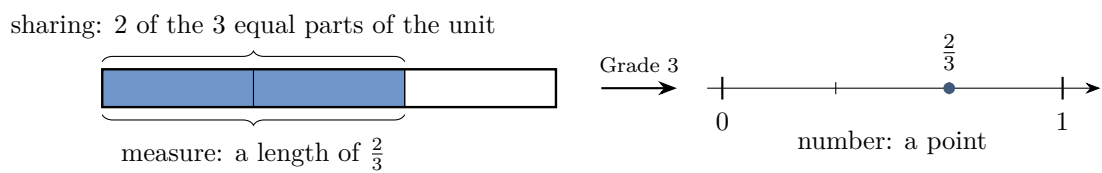

The bar is the object that belongs to both worlds (Figure~\ref{fig:bridge}). Cut into equal parts, it is a shared unit: ``2 of the 3 equal parts''. Laid along a line, its shaded part is a length: $\frac{2}{3}$, the bar itself, the unit, having length~$1$. The same figure supports both readings, and nothing has to be redrawn when the pupil moves from one to the other. Moreover, with the convention that the chapters follow for bars---a bar is cut across its width, into full-width strips---the two readings always coincide: every fair cut of a bar is then a length. This is specific to the bar cut in strips: a rectangle can be cut into three fair but L-shaped parts that are not lengths, and no fair cut of a disc, into sectors or otherwise, produces a part that is a length of the figure (Section~\ref{sec:disc}). Choosing the bar does not make the crossing happen---reading a share as a length remains something to be learnt---but it is the choice that makes both readings available on one figure.

The crossing is organised over three grades. In Grades~1 and~2 the reading is that of sharing: bars are cut and shaded, parts are counted and named, and the exercises also use discs and collections. In Grade~3 the lesson says it: ``the fraction is a length''. Pupils measure bars with rulers graduated in quarters, fifths or tenths of a unit, bound a length between two graduations, draw a segment of length $1+\frac{3}{10}$, and place fractions on the number line; lengths greater than one are already measured, written as a whole number of units and a fraction. In Grade~4 the single fraction greater than one, $\frac{5}{3}$, arrives: what changes is the notation and the iteration of the unit, not the possibility of a length above one. Evidence of the intended transition would be the pupil's ability to identify the reference unit and coordinate a partitioned bar, its numerical measure and the corresponding point when the representation changes. Counting unit fractions remains a legitimate strategy; its absence is not a criterion of understanding. That the two readings coexist on one figure, however, says nothing about how a pupil comes to coordinate them; the next subsection takes the crossing at the level of tasks.

\subsection{Crossing the bridge: the unit, task by task}\label{sec:crossing}

Reading a shaded bar as a measured length requires several distinct pieces of knowledge, which the bar does not supply by itself: that the fraction is a quantity and not a count of shaded pieces; that the unit is a \emph{specified} object---the strip, or the segment from $0$ to $1$, that counts as~$1$---and not the bar being measured; that the length is a displacement from~$0$; and that the endpoint is a number independent of the physical bar. The coordination of units that these require is a known difficulty \citep{steffe2010,hackenberg2007}, the forming and re-forming of units is central to the rational-number constructs \citep{behr1992}, and taking the wrong unit is the most frequent error on number-line tasks \citep{bright1988}. Two tasks of the chapters show how the progression addresses them; a third, which the literature on units would suggest, is excluded by design, and the reason for its exclusion is part of the bridge.

\paragraph{Task A: make the ruler, measure the bar (Grade~3).} The pupil receives a unit strip, folds it into four equal parts---``the strip becomes a ruler graduated in quarters''---lays its $0$ at one end of a bar and reads the graduation at the other end. The knowledge at stake is that the designated strip supplies the unit independently of the bar being measured. The chapter also includes lengths greater than one, such as $2+\frac{1}{4}$, which can be measured by iterating the reference unit. These lengths make a measurement interpretation useful, but do not by themselves force a pupil to abandon a sharing interpretation. Such a claim would require specifying how an inappropriate answer can be tested against the same reference.

\paragraph{Strategies and validation in Task~A.} Three strategies illustrate the distinction between a mathematical contradiction and feedback available to the pupil. (S1)~\emph{The measured bar taken as the unit.} A pupil may assign every complete bar measure~$1$, or describe a share of that bar without using the designated strip. The inconsistency becomes observable only when the reported measure is reconstructed using the designated strip and compared with the target. A bar longer than the strip does not perform this check by itself. (S2)~\emph{Marks instead of intervals.} Counting the mark at $0$ gives $\frac{3}{4}$ for a bar two quarters long \citep{bright1988}. Repeating the same mark-counting routine while drawing may reproduce the original length, so a read--draw check alone is insufficient. One possible check is to cover the target with quarter-unit pieces and count the pieces between the folds. (S3)~\emph{The entire drawn ruler taken as the unit.} On a ruler from $0$ to $2$ graduated in quarters, a bar five quarters long may be read as $\frac{5}{8}$ \citep{bright1988}. The mark~$1$ identifies the intended reference, but the pupil must use it to test the answer. In the Grade~3 reading exercise with a bar longer than the unit, the supplied answer format ``$1+\ldots$'' scaffolds this choice rather than requiring the pupil to detect the contradiction unaided. The online platform accepts or rejects an entered answer; this is external validation. The chapters do not ask for the reconstruction and covering checks described here. They are proposed additions that could provide feedback from the material situation, provided the teacher establishes the common reference and the pupil uses it. The anticipated knowledge is that a measure counts iterations of a specified unit or subunit; the effectiveness of these checks remains to be investigated.

\paragraph{Task B: one bar, two units (excluded from the initial progression).} A contrasting task would measure one unchanged bar with two reference strips: $\frac{3}{4}$ with the longer strip and $\frac{3}{2}$ with a strip half as long. Both results are measures. No conventional units and no explicit scalar-operator notation are needed to pose the problem. The progression deliberately postpones this coordination so that pupils can first work on magnitude, equivalence and operations within a stable reference, for example $\frac{1}{3}+\frac{1}{2}=\frac{5}{6}$. The intended benefit is to limit simultaneous conceptual demands, not an established reduction in cognitive load. Three changes must be distinguished: resizing an entire diagram preserves its represented fractions; changing the reference for an unchanged segment changes its numerical measure; applying $\frac{2}{3}$ to $18$~cm scales a quantity. Only the second is the task excluded here. In the linear model, the reference strip or the interval from $0$ to $1$ has measure~$1$; this does not identify a physical length with the abstract number~$1$.

\paragraph{Scope and cost of the stable reference.} The fraction chapters of Grades~1--3 do not specify lengths in conventional units such as centimetres. In Grades~4--5, examples such as $\frac{5}{10}$ of a metre, a lap of $\frac{1}{5}$~km and quarter-metre pieces from a $4$~m ribbon retain a named reference unit. Explicit applications such as ``$\frac{2}{3}$ of $10$~cm'' appear in Grade~6 alongside the written operator. This sequencing does not entail that operator instruction teaches changes of measurement unit, nor that placing fractions correctly guarantees understanding of such changes. The design postpones that learning and must evaluate it separately. Hypothesis~\ref{h:units} tests both whether numerical meaning transfers across drawings and whether pupils later coordinate a fixed length with different references.

\paragraph{Task C: draw a length with different graduations (Grade~3).} The pupil must draw a segment $\frac{1}{2}$ of a unit long with a ruler graduated in tenths. The method says: ``if the fraction is not written in the ruler's graduations, use an equivalent fraction''---$\frac{1}{2}$ is $\frac{5}{10}$, count five tenths. Here equivalence is not a rule on numerators and denominators but a necessity of measurement: the same length must be read in the ruler's parts. The anticipated error reads the denominator as a count and draws two tenths; a fold of the unit strip in half, laid against the drawing, would refute it. The same necessity drives the Grade~5 situation of two identical strips, one coloured over $5$ parts of $10$, the other over $1$ of $2$, with the question ``have they coloured the same length?''.

Tasks~A and~C, together with the placement task of the next subsection, address different parts of the intended transition: identifying the reference, expressing a length by a number, interpreting an endpoint and coordinating equivalent fractions. Their mathematical connection provides an organising principle. It does not show that pupils make the connection without instruction, or that the proposed material checks are already implemented. The initial tasks keep the designated reference stable; changes of reference remain a distinct later learning objective.

\subsection{From the bar to the point}

\begin{figure}[ht]
\centering
\begin{tikzpicture}
\begin{scope}[yshift=0.55cm]\fracbar{2}{3}{6}{0.35}\end{scope}
\numline{3}{6}{1}
\foreach \i in {1,2}{\draw[-{Stealth[length=1.6mm]},fracfill!60!black,thick] ({(\i-1)*2},0) to[out=-50,in=-130] node[below]{\scriptsize\i} ({\i*2},0);}
\fill[fracfill!60!black] (4,0) circle (2.2pt);
\node[above=2pt] at (4,0.95) {$\tfrac{2}{3}$};
\end{tikzpicture}
\hspace{0.8cm}
\begin{tikzpicture}
\begin{scope}[yshift=0.55cm]\fracbar{3}{2}{3}{0.35}\end{scope}
\numline{2}{3}{2}
\fill[fracfill!60!black] (4.5,0) circle (2.2pt);
\node[above=2pt] at (4.5,0.95) {$\tfrac{3}{2}$};
\end{tikzpicture}
\caption{From the bar to the point. Left: the bar $\frac{2}{3}$ above the number line; counting two thirds from $0$ reaches the point $\frac{2}{3}$. Right: the same routine continues past $1$, and $\frac{3}{2}$ (three halves) is a point between $1$ and $2$.}
\label{fig:line}
\end{figure}
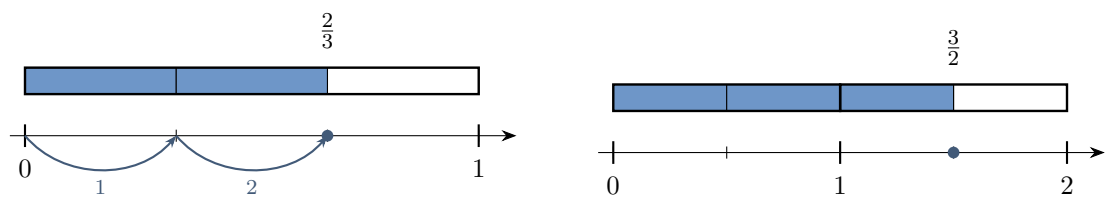

Placing the left end of a unit bar at $0$ on a number line whose unit is the length of that bar, the right end of its shaded part---the measured segment---falls on the point that is the fraction (Figure~\ref{fig:line}). In Grade~3 the placement is taught as the same three-step routine for every fraction: draw the line and mark $0$ and $1$; partition the segment from $0$ to $1$ into $b$ equal parts; count $a$ parts from $0$. Because the routine counts parts of a given size, it continues past $1$ once every unit segment, from $1$ to $2$ and beyond, is partitioned in the same way, and from Grade~4 single fractions greater than one occupy points beyond $1$ (Figure~\ref{fig:line}, right). The number line is thus the place where fractions and whole numbers meet on the same object, which is the central claim of the integrated theory \citep{siegler2011}.

\subsection{Comparison and equivalence}

\begin{figure}[ht]
\centering
\begin{tikzpicture}
\begin{scope}[yshift=0.8cm]\fracbar{2}{4}{4.8}{0.4}\node[left] at (-0.1,0.2){$\tfrac{2}{4}$};\end{scope}
\begin{scope}\fracbar{2}{3}{4.8}{0.4}\node[left] at (-0.1,0.2){$\tfrac{2}{3}$};\end{scope}
\draw[dashed] (2.4,1.3) -- (2.4,-0.2);\draw[dashed] (3.2,1.3) -- (3.2,-0.2);
\node at (2.4,-0.7) {(a) same numerator};
\end{tikzpicture}\hspace{0.8cm}
\begin{tikzpicture}
\begin{scope}[yshift=0.8cm]\fracbar{1}{3}{4.8}{0.4}\node[left] at (-0.1,0.2){$\tfrac{1}{3}$};\end{scope}
\begin{scope}\fracbar{2}{6}{4.8}{0.4}\node[left] at (-0.1,0.2){$\tfrac{2}{6}$};\end{scope}
\draw[dashed] (1.6,1.3) -- (1.6,-0.2);
\node at (2.4,-0.7) {(b) equivalence};
\end{tikzpicture}

\vspace{5mm}
\begin{tikzpicture}
\draw[-{Stealth[length=2mm]}] (-0.2,0) -- (8.5,0);
\foreach \i in {0,...,4}{\draw ({\i*2},0.1) -- ({\i*2},-0.1);}
\draw[thick] (0,0.15) -- (0,-0.15); \draw[thick] (8,0.15) -- (8,-0.15);
\node[above=3pt] at (2,0) {\small $\tfrac14$};
\node[above=3pt] at (4,0) {\small $\tfrac24$};
\node[above=3pt] at (6,0) {\small $\tfrac34$};
\node[below=3pt] at (0,-0.1) {\small $0$};
\node[below=3pt] at (4,-0.1) {\small $\tfrac12$};
\node[below=3pt] at (8,-0.1) {\small $1$};
\fill[fracfill!60!black] (4,0) circle (2.2pt);
\node at (4,-1.1) {(c) $\tfrac12$ and $\tfrac24$ mark the same point};
\end{tikzpicture}
\caption{Comparison and equivalence in the linear model. (a) With equal unit lengths and a common left end, $\frac{2}{4}<\frac{2}{3}$ is read as a difference of lengths. (b) $\frac{1}{3}$ and $\frac{2}{6}$ cover the same length. (c) Partitioning the unit into halves and into quarters shows that $\frac{1}{2}$ and $\frac{2}{4}$ are the same point.}
\label{fig:compare}
\end{figure}
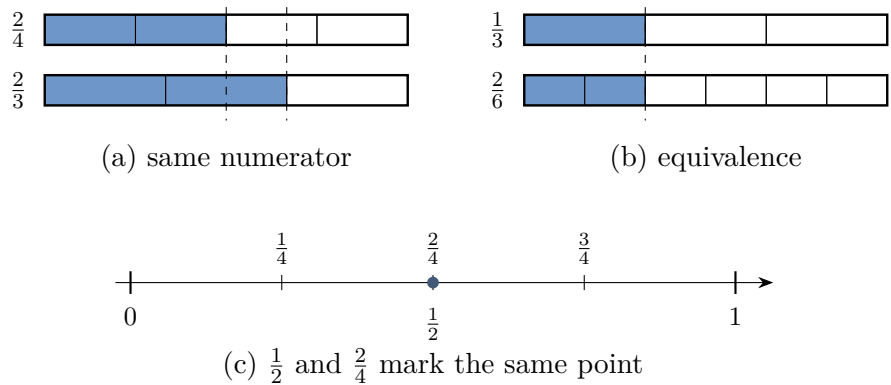

Two bars of the same unit length, aligned at their left end, compare by length (Figure~\ref{fig:compare}a). The comparison makes visible the fact behind the most common form of whole-number bias: with the same numerator, the larger denominator gives the shorter length, because each part is smaller. The progression introduces this in Grade~2 with unit fractions ($\frac{1}{6}<\frac{1}{3}$: the more parts, the smaller each part) and in Grade~3 on the number line, for fractions with the same denominator, with the same numerator ($\frac{5}{12}<\frac{5}{8}$), or with one denominator a multiple of the other; the general procedure through a common denominator follows in Grade~5, and benchmarks ($1$ and $\frac{1}{2}$) in Grade~6. Equivalence receives the same treatment: two bars partitioned differently can cover exactly the same length (Figure~\ref{fig:compare}b), and on the number line two fractions are equal when they mark the same point (Figure~\ref{fig:compare}c). Equality of fractions is thus an equality of magnitudes before it is a rule on numerators and denominators (Section~\ref{sec:p4}).

\subsection{Addition and subtraction}

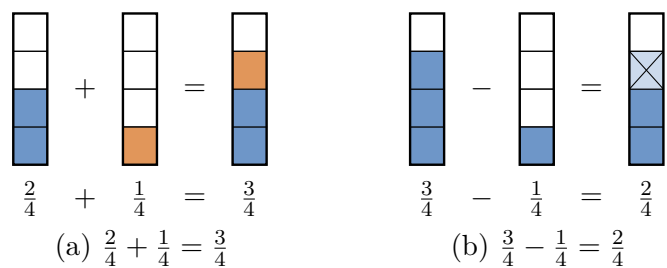
\begin{figure}[ht]
\centering
\begin{tikzpicture}[baseline=(current bounding box.center)]
\vbar{2}{0}{4}{fracfill}{white}
\node at (0.22,-0.45) {$\tfrac{2}{4}$};
\node at (0.95,1.0) {$+$}; \node at (0.95,-0.45) {$+$};
\begin{scope}[xshift=1.45cm]\vbar{1}{0}{4}{fracalt}{white}\node at (0.22,-0.45) {$\tfrac{1}{4}$};\end{scope}
\node at (2.4,1.0) {$=$}; \node at (2.4,-0.45) {$=$};
\begin{scope}[xshift=2.9cm]\vbar{2}{1}{4}{fracfill}{fracalt}\node at (0.22,-0.45) {$\tfrac{3}{4}$};\end{scope}
\node at (1.7,-1.1) {(a) $\tfrac24+\tfrac14=\tfrac34$};
\end{tikzpicture}
\hspace{1.6cm}
\begin{tikzpicture}[baseline=(current bounding box.center)]
\vbar{3}{0}{4}{fracfill}{white}
\node at (0.22,-0.45) {$\tfrac{3}{4}$};
\node at (0.95,1.0) {$-$}; \node at (0.95,-0.45) {$-$};
\begin{scope}[xshift=1.45cm]\vbar{1}{0}{4}{fracfill}{white}\node at (0.22,-0.45) {$\tfrac{1}{4}$};\end{scope}
\node at (2.4,1.0) {$=$}; \node at (2.4,-0.45) {$=$};
\begin{scope}[xshift=2.9cm]\vbarsub{3}{1}{4}{fracfill}\node at (0.22,-0.45) {$\tfrac{2}{4}$};\end{scope}
\node at (1.7,-1.1) {(b) $\tfrac34-\tfrac14=\tfrac24$};
\end{tikzpicture}
\caption{Addition and subtraction with the same denominator, in the layout of the lessons (Grade~3). The unit bar is cut into quarters; the parts of the two terms keep their colour in the result. (a) Parts of the same size are counted together. (b) The removed part is crossed out.}
\label{fig:operations}
\end{figure}

Adding lengths means placing them end to end. In the bar model the parts of $\frac{2}{4}$ and of $\frac{1}{4}$ have the same size, so the sum is counted in quarters, and the denominator does not change because the size of the parts does not change (Figure~\ref{fig:operations}a); subtraction removes parts of that size (Figure~\ref{fig:operations}b). When the parts differ in size, as in $\frac{1}{2}+\frac{1}{4}$, they cannot be counted together; one of the fractions is first re-partitioned so that all parts have the same size (Figure~\ref{fig:chain}). The search for a common denominator thus appears as a necessity of the model before it becomes an algorithm, and the classic error $\frac{1}{2}+\frac{1}{3}=\frac{2}{5}$ produces a length that visibly does not match the two pieces placed end to end. In the progression, sums with the same denominator appear in Grade~2 exercises. Grade~3 defines them and, for unlike denominators, states that the parts must first be made the same size and that the numerators cannot be added directly, and gives the general procedure through a common multiple of the denominators, with examples in which one denominator is a multiple of the other; Grade~5 applies the same procedure to any denominators. The written form of these calculations is the object of Principle~\ref{pr:chain} (Section~\ref{sec:p5}).

\subsection{Multiplication and division}

\begin{figure}[ht]
\centering
\begin{tikzpicture}[scale=1.05]
\foreach \i in {0,1,2}{\foreach \j in {0,1}{
  \pgfmathtruncatemacro{\filled}{and(\i<2,\j<1)}
  \ifnum\filled=1 \draw[partbox,fill=fracmix!45] (\i,\j*1.5) rectangle ++(1,1.5); \else \draw[partbox] (\i,\j*1.5) rectangle ++(1,1.5); \fi}}
\draw[unitbox] (0,0) rectangle (3,3);
\begin{scope}[yshift=-0.45cm]\fracbar{2}{3}{3}{0.25}\end{scope}
\begin{scope}[xshift=-0.2cm,rotate=90]\fracbar[shadedalt]{1}{2}{3}{0.25}\end{scope}
\node[below] at (1.5,-0.5) {$\tfrac23$};
\node[left] at (-0.5,0.75) {$\tfrac12$};
\node[right] at (3.3,1.5) {$\tfrac23\times\tfrac12=\tfrac{2\times1}{3\times2}=\tfrac26$};
\end{tikzpicture}
\caption{The area model of multiplication. The unit square is partitioned into $3\times 2$ equal cells; the rectangle with sides $\frac{2}{3}$ and $\frac{1}{2}$ is made of $2\times 1$ of them.}
\label{fig:area}
\end{figure}
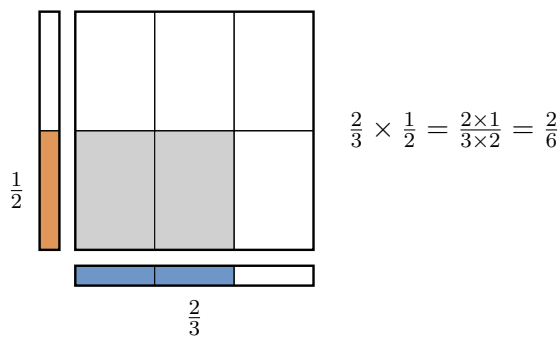

Multiplication is where the purely linear model reaches its limit, and the progression introduces three models in Grade~5 before the general rule of Grade~7. A whole number times a fraction is repeated addition on the number line ($3\times\frac{2}{3}=2$). A fraction of a fraction is shown on an area: half of one third of a square is one sixth of it, ``also written $\frac{1}{2}\times\frac{1}{3}$''. And dividing a whole number by a unit fraction is measured on the line: there are six halves in $3$, so $3\div\frac{1}{2}=6$. A short remark asks pupils to predict, before computing, whether the result will be larger or smaller than the starting number. In Grade~7 the general rule $\frac{a}{b}\times\frac{c}{d}=\frac{ac}{bd}$ is justified by an area (Figure~\ref{fig:area}): cutting the unit square into $b$ equal columns and $d$ equal rows gives $bd$ equal cells, each of area $\frac{1}{bd}$, and a rectangle with sides $\frac{a}{b}$ and $\frac{c}{d}$ is tiled by $a\times c$ such cells; when a factor is greater than one the rectangle extends beyond the unit square, but the tiling argument is unchanged. It covers non-negative fractions; signs are handled separately, by the sign rules of Grade~7. Mathematically the area model extends the linear backbone, since each side of the rectangle is a fraction bar; whether a pupil recognises in that side the bar he knows, rather than a new object, is not given by the continuity and would have to be observed; it does, however, introduce a change of dimension, which in Vergnaud's terms is a different multiplicative structure, the \emph{product of measures} \citep{vergnaud1983}. It also presupposes that the area of a rectangle with fractional sides is the product of its sides; in the progression this is prepared by the Grade~5 reading of ``of'' as multiplication, so that the area model illustrates the rule as much as it proves it, and the Grade~7 chapter's own first example with a factor greater than one is computed by the rule alone.

\subsection{Why the disc does not cross the bridge}\label{sec:disc}

\begin{figure}[ht]
\centering
\begin{tikzpicture}
\begin{scope}
\fill[fracfill] (0,0) -- (0:1) arc (0:120:1) -- cycle;
\fill[fracfill] (0,0) -- (120:1) arc (120:240:1) -- cycle;
\draw[line width=0.9pt] (0,0) circle (1);
\foreach \a in {0,120,240}{\draw[line width=0.4pt] (0,0) -- (\a:1);}
\node at (0,-1.6) {\small a share, not a length};
\end{scope}
\begin{scope}[xshift=4cm,yshift=-0.25cm]
\fracbar{2}{3}{4.5}{0.5}
\draw[decorate,decoration={brace,amplitude=4pt,mirror}] (0,-0.1) -- (3,-0.1) node[midway,below=4pt]{\small a length of $\tfrac23$};
\node at (2.25,-1.35) {\small a share and a length};
\end{scope}
\end{tikzpicture}
\caption{Two of three equal parts, on a disc and on a bar. Both figures show a fair share; only the bar shows a length that can be placed on a number line.}
\label{fig:disc}
\end{figure}
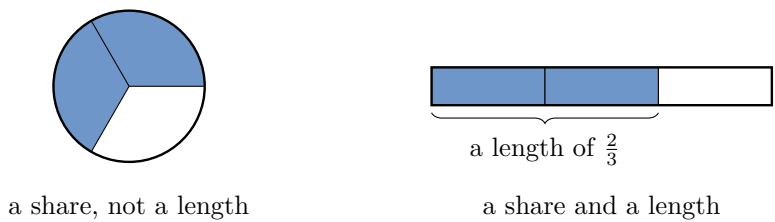

Circular area models (``pizzas'') are one of the representations that textbooks use for fractions, in proportions that vary between countries \citep{alajmi2012}, and they are present in the progression: the first problem of Grade~1 is a pizza, halves and quarters are shown on discs as well as on rectangles and strips, and the exercises of Grades~1 to~4 ask pupils to recognise fractions in many kinds of figures, discs included. A disc is a good object for sharing: a cake is round, and fair cuts are easy to judge. But a sector is a share and not a length (Figure~\ref{fig:disc}); \citet{wu2014} makes the same point when he writes that a pizza is not a good model for the whole once fractions are to be numbers. It does not lie on an axis with the whole numbers, so the passage to the number line requires a translation that the bar does not; fractions greater than one require several discs, with no axis on which to continue counting; comparing sectors means comparing areas or arcs rather than lengths from a common origin; and although a fraction of a fraction can be shown on a sector by successive partitions, the product of two fractions as the area of a rectangle with fractional sides has no direct circular counterpart. A disc can represent a magnitude---an area---and sectors can be compared; but it does not support the linear transition as directly as the bar, and that, rather than any defect as a model of sharing, is why the progression does not let it carry the definition, the number line or the operations.

These are arguments from the design, not evidence that circles fail. The Rational Number Project, whose curriculum outperformed commercial curricula in a large comparison \citep{cramer2002}, found fraction circles to be the most powerful of its models for building mental images of fraction size \citep{cramer1995,cramer2002b}. Experimental comparisons that favour the number line over area models for magnitude understanding \citep{hamdan2017,gunderson2019} support the choice of a linear backbone, but they do not show that circles should be absent. The position taken here is therefore one of hierarchy, not exclusion: the bar carries the meaning and the operations, and other models are related back to it.

% ======================================================================
\section{Principle 3: establishing numbers before explicit operators}\label{sec:p3}
% ======================================================================

\begin{principle}[Numbers before explicit operators]\label{pr:ratio}
A fraction is first taught as a number. The progression postpones the explicit conceptualisation and the symbolic treatment of the fraction as a scalar operator---the scaling $x\mapsto\frac{a}{b}\,x$, written $\frac{a}{b}\times N$---not the application of fractions to quantities, which fair sharing provides from Grade~1. ``Two thirds of thirty'' is treated as a situation, a part--whole relation on a quantity in which the number $\frac{2}{3}$ is applied. Such situations, and part--part ratios written $a:b$, are developed in a companion \emph{ratio} strand, where they receive their meaning from fair sharing; the fraction strand prioritises numerical meaning and schedules explicit operator notation after work on magnitude, equivalence and arithmetic. This order is intended to consolidate numerical meaning; the schedule alone does not establish that pupils have achieved it. The order is one of institutionalisation, not of encounter: a meaning may be met, or used locally, in either strand before it is stated (Section~\ref{sec:encounter}).
\end{principle}

\subsection{The operator as a situation}

The operator is one of the recognised sub-constructs of the rational number (Section~\ref{sec:background}), and a rational number $r$ and the map $x\mapsto rx$ are two faces of one object; the design choice is therefore one of \emph{sequencing}, not of nature. What is postponed is not the application of a fraction to a quantity---half of twelve strawberries is found by sharing in Grade~1, a third of twelve beads in Grade~2---but its conceptualisation as the action of a number and its writing as a multiplication: in the early tasks the pupil shares a collection and counts a group, and no fraction symbol acts on anything; from Grade~5 in the ratio strand and Grade~6 in the fraction strand, the same relation is written $\frac{a}{b}\times N$ and the number is seen to act. The distinction is thus didactic and symbolic, not mathematical. Many presentations give the fraction the status of an operator---``$\frac{2}{3}$ of''---alongside, or before, its status as a number. The risk, as the design sees it, is that pupils know what $\frac{2}{3}$ \emph{does} before they know what it \emph{is}: ``$\frac{2}{3}$ of the whole'' then names an action on a referent, and the question ``where is $\frac{2}{3}$?'' has no answer for them. The progression takes the opposite order. The intended order establishes numerical meaning first, through the bridge of Principle~\ref{pr:linear}; ``of'' applied to a quantity is treated as a situation, a relation between a part and a whole, which belongs to the ratio strand; and within the fraction strand the word ``of'' first appears between numbers, in Grade~5, where ``half of one third'' is used to model the product of two fractions (Section~\ref{sec:p2}), not to act on a quantity.

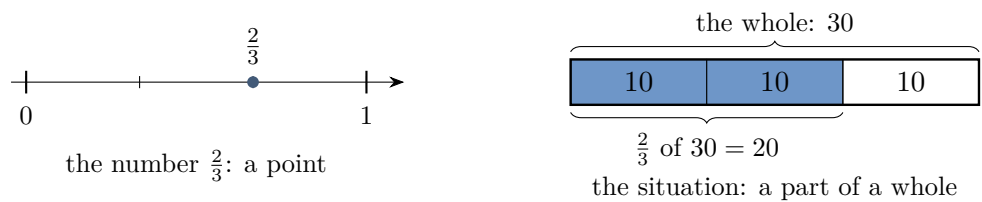
\begin{figure}[ht]
\centering
\begin{tikzpicture}
\begin{scope}
\numline{3}{4.5}{1}
\fill[fracfill!60!black] (3,0) circle (2.2pt);
\node[above=3pt] at (3,0) {$\tfrac{2}{3}$};
\node at (2.25,-1.1) {\small the number $\tfrac23$: a point};
\end{scope}
\begin{scope}[xshift=7.2cm,yshift=-0.3cm]
\foreach \i in {0,1}{\draw[partbox,fill=fracfill] (\i*1.8,0) rectangle ++(1.8,0.6) node[midway]{$10$};}
\draw[partbox] (3.6,0) rectangle ++(1.8,0.6) node[midway]{$10$};
\draw[unitbox] (0,0) rectangle (5.4,0.6);
\draw[decorate,decoration={brace,amplitude=4pt}] (0,0.7) -- (5.4,0.7) node[midway,above=4pt]{\small the whole: $30$};
\draw[decorate,decoration={brace,amplitude=4pt,mirror}] (0,-0.1) -- (3.6,-0.1) node[midway,below=4pt]{\small $\tfrac23$ of $30 = 20$};
\node at (2.7,-1.1) {\small the situation: a part of a whole};
\end{scope}
\end{tikzpicture}
\caption{The number and the situation. Left: $\frac{2}{3}$ is a point on the number line. Right: ``$\frac{2}{3}$ of $30$'' is a part--whole situation on a quantity: the whole is shared into three equal parts and two are taken. The design introduces the second as an application of the number represented on the left.}
\label{fig:operator}
\end{figure}

That strand introduces the situation in Grade~3, under the heading \emph{Part in part--whole relationships}, with a collection as the whole: ``to find the number of apples corresponding to $\frac{1}{2}$ of 4 apples, divide the whole into 2 equal parts and select 1 part; count the apples in the selected part.'' The meaning is fair sharing of a whole, the same action as in Grade~1, now applied to a quantity. In Grade~5 the ratio strand states the link with multiplication: ``half of'' can be read as a sharing ($4\div 2=2$) or as a multiplication ($\frac{1}{2}\times 4=2$), and ``the sharing reading is the one that computes; the multiplication reading is the one that writes.'' The fraction strand itself writes ``of'' as a multiplication between two fractions in Grade~5 (``half of one third \ldots{} also written $\frac{1}{2}\times\frac{1}{3}$''), to model the product; it writes the operator on a quantity only in Grade~6: ``$\frac{2}{3}$ of 18---split 18 into 3 equal parts, one part is $18\div 3=6$, two parts make $2\times 6=12$; taking $\frac{a}{b}$ of a number always works this way: divide by the denominator, multiply by the numerator, which is exactly the multiplication $\frac{a}{b}\times N$,'' followed by the proposition $\frac{a}{b}\times N=\frac{aN}{b}$ and its two readings, divide first or multiply first. By then the pupil has placed $\frac{2}{3}$ on a line (Grade~3), gone beyond one (Grade~4) and added and compared fractions (Grade~5); the operator is intended to arrive as a use of a number already studied (Figure~\ref{fig:operator}). Whether its numerical status is secure requires assessment, not inference from the list of topics taught.

\subsection{Part--part ratios and the mediant}

\begin{figure}[ht]
\centering
\begin{tikzpicture}
\foreach \i in {0,1,2}{\draw[partbox,fill=fracalt] (\i*1.2,0) rectangle ++(1.2,0.6);}
\foreach \i in {3,4}{\draw[partbox,fill=fracfill!45] (\i*1.2,0) rectangle ++(1.2,0.6);}
\draw[unitbox] (0,0) rectangle (6,0.6);
\draw[decorate,decoration={brace,amplitude=4pt}] (0,0.7) -- (3.6,0.7) node[midway,above=4pt]{\small cherries: 3 parts};
\draw[decorate,decoration={brace,amplitude=4pt}] (3.6,0.7) -- (6,0.7) node[midway,above=4pt]{\small apples: 2 parts};
\draw[decorate,decoration={brace,amplitude=4pt,mirror}] (0,-0.1) -- (6,-0.1) node[midway,below=4pt]{\small all the fruit: 5 parts};
\node[right,align=left,font=\small] at (6.4,0.3) {part--part: cherries $:$ apples $=3:2$,\\ \hphantom{part--part: }cherries $=\tfrac32$ of the apples\\ part--whole: cherries $=\tfrac35$ of the fruit};
\end{tikzpicture}
\caption{Part--part and part--whole readings of the same situation, as distinguished in the ratio strand.}
\label{fig:ratio}
\end{figure}
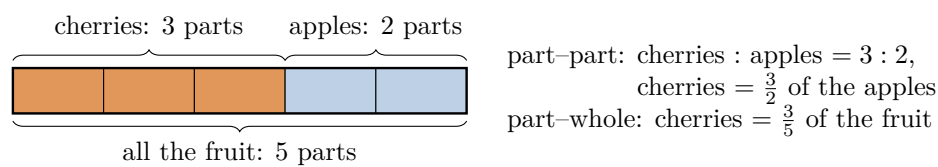

The ratio strand also carries the comparison of two parts. A ratio is defined in its Grade~2 as a comparison of two quantities, written $a:b$; equivalent ratios follow, then the distinction between part--part and part--whole ratios in Grade~4, proportion in Grade~5, and the unitary and cross-multiplication methods in Grade~6. The Grade~4 chapter states what the fraction measures in each case (Figure~\ref{fig:ratio}). In the part--part ratio cherries$\,:\,$apples $=3:2$, the fraction $\frac{3}{2}$ is the quotient cherries $\div$ apples---there are $\frac{3}{2}$ as many cherries as apples---and it is not a part of the whole: the cherries are $\frac{3}{3+2}=\frac{3}{5}$ of the fruit. Reading a part--part ratio $a:b$ as ``$\frac{a}{b}$ of the total'' is the error that this explicit statement targets. The remark thus meets the quotient reading two years before the fraction strand institutionalises the quotient (Grade~6); Section~\ref{sec:encounter} explains how the model treats such encounters.

The \emph{mediant} illustrates why the numerical value of a fraction must be distinguished from the quantities in a ratio situation; it does not by itself determine their instructional order. When two ratios describe collections that are put together, with the same referents, adding numerators and denominators gives the ratio of the combined collection: one red marble out of two, put together with one red marble out of three, gives two red marbles out of five; a batch of one cup of syrup for two of water, mixed with a batch of two cups of syrup for three of water, gives three cups of syrup for five of water. The mediant is an operation on the actual quantities, not on the ratios as classes: pooling a batch of two cups of syrup and three of water with a batch containing two cups of syrup and four of water gives $4:7$; pooling it with one cup of syrup and two of water gives $3:5$. The ratios $2:4$ and $1:2$ are equivalent, but the actual batch sizes differ. For fractions as numbers, $\frac{1}{2}+\frac{1}{3}=\frac{2}{5}$ is false. The same written operation is right in one role and wrong in the other, which is why a design that writes ratios as fractions has to say, each time, which role is meant; teachers and textbooks model the relation between ratios and fractions in several incompatible ways \citep{clark2003}.

In Vergnaud's analysis of the isomorphism of measures \citep{vergnaud1983}, two measure spaces are related proportionally---cakes and eggs, for instance, with four eggs per cake. A \emph{scalar} relation links two values within one space and is dimensionless: two cakes are twice one cake, so eight eggs are twice four. A \emph{function} relation links the two spaces and carries a quotient of units: four eggs per cake. The operator ``$\frac{2}{3}$ of $30$'' is a scalar operator; a rate such as four eggs per cake is a function relation; a part--part ratio such as cherries$\,:\,$apples $=3:2$ is a multiplicative comparison between two cardinalities, which, read as ``$\frac{3}{2}$ cherry per apple'', takes the form of a function relation between two measure spaces whose units differ, even though both quantities are counts. It should be said that Vergnaud's own framing runs the other way: for him fractions, even unit fractions, are relationships between quantities rather than quantities \citep{vergnaud1983}. Principle~\ref{pr:ratio} is a didactic choice against that framing, not a reading of it: it hypothesises that explicit coordination of these relations benefits from prior work on numerical magnitude. Separate strands make this intention visible in the document architecture, while the ratio strand can begin with whole-number comparisons. Their separation does not establish the order of pupils' encounters, which must be considered across both strands.

\subsection{Meeting a meaning, institutionalising it}\label{sec:encounter}

A pupil who follows both strands meets several meanings of the fraction notation before the fraction strand states them: separating documents does not separate pupils' experiences. Three relations to a meaning can be distinguished for the documentary analysis: \emph{meeting} it in a situation without an explicit general statement; \emph{using} it locally to solve a particular problem; and \emph{presenting it as reference knowledge} through a definition, proposition or general method. The last is an intended opportunity for institutionalisation in the sense of the theory of didactical situations \citep{brousseau1997}, not evidence that classroom institutionalisation has occurred. The distinction is related to \citet{douady1986}'s account of a notion used as a tool and studied as an object. Principle~\ref{pr:ratio} proposes an order for explicit reference knowledge: numerical magnitude on the line from Grade~3, explicit multiplication by a fraction in the ratio strand's Grade~5 and the fraction strand's Grade~6, and the quotient characterisation in Grade~6. Earlier sharing encounters remain part of the design.

This distinction does not automatically remove sequencing discrepancies. The Grade~4 ratio remark reads $\frac{3}{2}$ as cherries $\div$ apples; its position in a remark supports classifying it as a local use, but its instructional role depends on whether later explanations and tasks require pupils to use that reading independently. Likewise, the Grade~3 method for finding a part of a collection states a sharing procedure, but can contribute to an operator conception even without multiplication notation. The Grade~2 ratio definition explicitly writes $a:b$ ``or as the fraction $\frac{a}{b}$'', introducing a general link before the source progression's number-line work and part--part/part--whole distinction (Section~\ref{sec:discussion}, item~\ref{it:ratiostrand}). The intended sequencing therefore requires coordination between strands and an explicit return to earlier encounters, including the Grade~4 quotient remark. Block types are documentary indicators, not sufficient evidence of how a meaning functions for a pupil. Whether early encounters help or hinder later consolidation remains an empirical question.

\subsection{The French programmes on this point}

The French programmes make a different choice, and it is worth stating precisely, because the difference is smaller than a first reading suggests. The cycle~3 programme in force since 2025 lists, among the meanings that fractions take in CM1, that ``fractions also acquire the status of multiplicative operator, for the particular case of unit fractions; pupils learn to compute fractions of quantities or magnitudes, such as a third of 12 marbles'', extended in CM2 to any fraction (``two thirds of 12~€''); in \emph{sixi\`eme} the objective is ``to use a multiplication to apply a fraction to a whole number'', the programme noting that to make a fraction operate ``not only on a quantity or a magnitude, as in CM, but also on a whole number \ldots{} constitutes a higher level of abstraction'' \citep{men2025cycle3}. The programme thus grades the operator as the progression does---quantities first, numbers later---and the accompanying guidance gives it its meaning through sharing (``to find a third of twelve eggs, I share into three equal parts'') \citep{men2025fractions}. The term \emph{ratio} appears only in the cycle~4 programme of 2026, with proportionality; cycle~3 names the part--whole relation a \emph{rapport} from CM1 on, and cycle~4 defines a fraction as ``the quotient of two integers, which can be seen as a number, an expression, and also as an operator (fraction of a quantity)'' \citep{men2026cycle4}.

The difference lies in two places. First, in CM1 the operator is one of the meanings that the fraction itself acquires, in the same year as the half-line that gives it the status of a number; in the progression the operator on a quantity is conceptualised in the ratio strand, and the fraction strand writes it as a multiplication in Grade~6, after several source levels devoted to numerical meaning. Second, in the programmes ratio and operator are \emph{meanings} of the fraction; here they are \emph{situations} in which one number is used. Our design hypothesis is that the second order gives ``$\frac{2}{3}$ of 30'' its meaning as a use of a number already placed on the line, whereas a fraction that is at once a point, a part and an operator asks the pupil to hold three statuses in the same year; whether this proposed sequencing matters for learning is Hypothesis~\ref{h:operator}. Presenting ratio and operator as legitimate readings in a later grade is compatible with that aim; the audit must examine when and how the readings are introduced and connected, rather than treating their vocabulary alone as a contradiction (Section~\ref{sec:discussion}).

% ======================================================================
\section{Principle 4: technique before formalism}\label{sec:p4}
% ======================================================================

\begin{principle}[Technique before formalism]\label{pr:technique}
A result needed before its full theoretical status is accessible should receive an informal explanation in the model and a mathematically valid local justification using available knowledge; its formal generalisation may be deferred. Validity conditions must be respected from the first use and made explicit in language appropriate to the task. An activity can embody a condition before pupils formulate it verbally. The audit distinguishes this implicit constraint from an explicit statement and identifies where a local justification is still a proposed revision.
\end{principle}

\subsection{Equivalence: a rule justified in the model}

Equivalence is introduced in Grade~3, in the same year as the number line, with its justification in the model. Two fractions are equivalent when they represent the same amount; cutting every part into the same number of smaller parts (a number of parts, hence at least $1$---our gloss) multiplies the numerator and the denominator by that number, and grouping parts back divides both by a number that must divide them exactly. Pupils also recognise equivalent fractions as the same point on the number line (Figure~\ref{fig:compare}c). The rule is thus a description of an operation on the partition, not an arbitrary manipulation of two numbers; it is restated in Grade~5 with the condition that the multiplier be nonzero, and in Grade~6 equal fractions are presented as two ways of writing the same number.

\subsection{Cross-multiplication: a criterion, not a construction}

In Grade~6 the progression states the criterion
\[
  \frac{a}{b}=\frac{c}{d}\quad\text{exactly when}\quad a\times d=b\times c ,
\]
with no explicit hypotheses beside the display (the denominators being positive integers under the definition then in force). Its remark multiplies both sides by $b\times d$ and says that the same reasoning works in reverse. The identity $\frac{a}{b}\times b=a$ need not depend on the later quotient interpretation: regrouping $ab$ pieces of size $\frac{1}{b}$ gives $a$ units. Similarly, multiplication by a positive integer preserves and reflects equality of lengths. Thus the explanation is not inherently circular. The documentary issue is that these justifications are not developed in the remark, whose cancellation language anticipates the quotient treatment later in the chapter. A proof already supported by the earlier equivalence rule is available.

\paragraph{Recommended local proof by subdivision.} For non-negative integers $a,c$ and positive integers $b,d$, subdividing each $b$-th into $d$ equal pieces and each $d$-th into $b$ equal pieces gives
\[
  \frac{a}{b}=\frac{ad}{bd},\qquad
  \frac{c}{d}=\frac{bc}{bd}.
\]
Both quantities now count parts of the same size $\frac{1}{bd}$. They are equal exactly when their counts are equal, hence
\[
  \frac{a}{b}=\frac{c}{d}\quad\Longleftrightarrow\quad ad=bc.
\]
This proves both directions without changing the unit or assuming the later quotient definition. It is a proposed replacement for the current Grade~6 remark, not a proof already present in that chapter. The criterion is used to test equality and solve missing-term equations such as $\frac{10}{5}=\frac{x}{8}$. Grade~7 restates it with $b,d\neq0$ and proves the direct implication; Grade~8 states it for algebraic expressions on their permitted domains.

In the praxeological terms of \citet{chevallard1999}, the Grade~6 chapter offers a type of task (decide whether two fractions are equal; find a missing term), a technique (compare the cross products) and a technology (the remark on why it works), while the theory that would justify the technology is left implicit. Two such theories are available. The older is the theory of proportion: for positive whole numbers, Euclid's \emph{Elements} prove that four numbers are proportional exactly when the product of the first and fourth equals the product of the second and third (Book~VII, Proposition~19; \citealp{euclidheath1926}). The more recent is the construction of $\mathbb{Q}$ as the field of fractions of $\mathbb{Z}$, in which the relation $(a,b)\sim(c,d)\iff ad=bc$ on $\mathbb{Z}\times(\mathbb{Z}\setminus\{0\})$ is the \emph{definition} of equality and the rational numbers are its equivalence classes. There, the criterion is not a consequence of a prior meaning of $\frac{a}{b}$; it is what gives the symbol a meaning. In the progression the order is reversed: the fraction already has a meaning as a magnitude, and the criterion is a property of magnitudes.

In terms of didactic transposition \citep{chevallard1985}, the source chapter offers a technique and a brief local explanation; the subdivision proof would supply the explicit valid justification required by Principle~\ref{pr:technique}. Deferring the construction of $\mathbb{Q}$ as a quotient set is a legitimate decision about formalisation. Leaving an available two-direction argument unstated is a separate implementation gap. The proposed revision concerns that gap and does not require teaching the field construction at this level.

\subsection{The quotient: a characterisation stated from one example}

The hand-over of Section~\ref{sec:p1} follows the same pattern. Grade~6 introduces $a\div b=\frac{a}{b}$ with one sharing situation (two cakes shared among three people), states it for whole-number numerators and positive whole-number denominators, and infers that $\frac{a}{b}$ is the number which, multiplied by $b$, gives $a$. A general argument is available within the progression---for $a\geq0$ and $b\geq1$, sharing $a$ units among $b$ sharers gives each sharer $a$ pieces of size $\frac{1}{b}$, that is $a\times\frac{1}{b}=\frac{a}{b}$, a product written in the same chapter---but the chapter does not give it, so that the statement on which the later grades rest functions as a technique whose technology is a single example. A progression that applies Principle~\ref{pr:technique} consistently would make that argument explicit at the hand-over, since it is the point where the operative definition changes.

\subsection{Conditions made explicit when they start to matter}

Conditions of validity determine whether an expression is defined. In Grade~1, cutting a unit into equal parts embodies a positive-integer denominator, but the chapter does not formulate the condition. In Grade~4 it is stated in words: ``the denominator is never zero: a unit cannot be divided into zero equal parts''. Grade~5 attaches a nonzero multiplier to the equivalence rule, and Grade~7 writes $b\neq0$ in the definition of rational numbers and as a hypothesis of most rules. These are three distinct forms of treatment: a constraint inherent in an activity, a verbal statement and a symbolic condition. They should not be described as equally explicit from the outset. In Grade~8, nonzero denominators depend on variable values and become an object of study. The identities $\frac{9x^3}{30x}=\frac{3x^2}{10}$ for $x\neq0$ and $\frac{x^2-9}{x-3}=x+3$ for $x\neq3$ preserve values on the original domains; the expressions on the two sides have different natural domains. The lessons write excluded values beside cancellations. An answer that inherits an explicit domain from its question need not be false if it does not repeat it, but restating the domain makes the intended habit visible and avoids an unrestricted reading.

% ======================================================================
\section{Principle 5: one written form, the chain of equalities}\label{sec:p5}
% ======================================================================

\begin{principle}[The chain of equalities]\label{pr:chain}
The intended convention for a multi-step calculation is a chain of equal expressions, with one transformation per step. Chains appear in Grades~2--3; written reasons are introduced in lessons from Grade~4 and occur unevenly in model answers. The vertical layout is intended to support later algebraic work, while equality of expressions and equivalence of equations require separate teaching.
\end{principle}

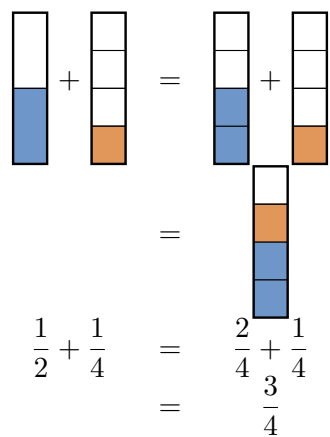
\begin{figure}[ht]
\centering
\renewcommand{\arraystretch}{1.3}
\begin{tabular}{ccc}
\begin{tikzpicture}[baseline=(current bounding box.center)]\vbar{1}{0}{2}{fracfill}{white}\end{tikzpicture} $+$
\begin{tikzpicture}[baseline=(current bounding box.center)]\vbar{1}{0}{4}{fracalt}{white}\end{tikzpicture}
& $=$ &
\begin{tikzpicture}[baseline=(current bounding box.center)]\vbar{2}{0}{4}{fracfill}{white}\end{tikzpicture} $+$
\begin{tikzpicture}[baseline=(current bounding box.center)]\vbar{1}{0}{4}{fracalt}{white}\end{tikzpicture} \\
& $=$ &
\begin{tikzpicture}[baseline=(current bounding box.center)]\vbar{2}{1}{4}{fracfill}{fracalt}\end{tikzpicture} \\[4pt]
$\dfrac{1}{2}+\dfrac{1}{4}$ & $=$ & $\dfrac{2}{4}+\dfrac{1}{4}$ \\
 & $=$ & $\dfrac{3}{4}$ \\
\end{tabular}
\caption{The Grade~3 lesson's figure for $\frac{1}{2}+\frac{1}{4}$: the bars above, the chain of equalities below, line for line. The half is first re-cut into quarters, then the quarters are counted; each written line records one of the two pictures.}
\label{fig:chain}
\end{figure}

Fractions are among the first places in school arithmetic where a written calculation routinely has steps of different kinds that cannot be set out in columns: to add $\frac{1}{2}$ and $\frac{1}{4}$ one must first rewrite $\frac{1}{2}$ as $\frac{2}{4}$, then count the quarters. The progression uses this moment to install a written form that can also support algebra. The first chains appear in the answers to Grade~2 exercises, $\frac{3}{8}+\frac{4}{8}=\frac{3+4}{8}=\frac{7}{8}$, each written next to its sentence (``three eighths plus four eighths make seven eighths''). In the Grade~3 lesson the chain is set under the bars that it describes, line for line (Figure~\ref{fig:chain}): the picture gives each line its meaning, and the line records what the picture shows. From Grade~4 the lessons write a reason in the margin, here the one for the decomposition:
\[
\begin{aligned}[t]
\frac{9}{4} &= \frac{8}{4}+\frac{1}{4} && \text{(8 quarters make 2 units)}\\
            &= 2+\frac{1}{4}
\end{aligned}
\]
and the same layout serves, unchanged, the simplification of $\frac{-4}{-6}$ in Grade~7, each line annotated, and the algebraic simplifications of Grade~8, where the condition that accompanies them (``valid for $x\neq 0$'' for $\frac{9x^3}{30x}$) is of a new kind: not a rule applied but a restriction on the letter. The chapters carry the principle unevenly, however: the lessons annotate most of their chains from Grade~4, but the answers of the exercise sets annotate only a minority of theirs before Grade~8 ($4$ of the $48$ chains of Grade~4, $15$ of $106$ in Grade~6, $44$ of $164$ in Grade~7, against $69$ of $93$ in Grade~8; Table~\ref{tab:trace}), so that in the middle grades the model answers do not model the habit; whether the exercises ask pupils for reasons cannot be read from the answers alone (Section~\ref{sec:discussion}).

The rationale is structural. An equation is solved by a chain of equivalent equations, a proof by a chain of justified statements, an algebraic simplification by a chain of equal expressions. The three share one layout: a vertical sequence in which each line follows from the previous one by a stated reason, so that the reader can check the step rather than the result. The relation between lines is not the same in the three cases---equality of expressions, equivalence of equations under conditions on the unknown, implication between statements. In
\[
\frac{1}{2}+\frac{1}{4}=\frac{2}{4}+\frac{1}{4}=\frac{3}{4}
\qquad\text{and}\qquad
\frac{x}{2}=3\iff x=6 ,
\]
the layout is the same and the sign is not: the first chain joins expressions that name one number, the second joins equations that have the same solutions, and confusing the relations can produce a string such as $\frac{x}{2}=3=x=6$. The Grade~8 domain condition is where the difference first becomes visible in the progression. What the design aims to carry over from Grade~3, with reasons developed from Grade~4, is the layout and attention to each transformation; the logic of each relation must be taught; that the shared layout helps them learn the second relation is a mathematical continuity seen by the expert, and a hypothesis about the learner (Hypothesis~\ref{h:chain}). Asking for reasons may make transformations available for discussion, including the invalid step $\frac{1}{2}+\frac{1}{3}=\frac{2}{5}$ and omitted domain conditions. However, a pupil can name an invalid rule, and a correct single-line calculation can be justified. The convention therefore needs mathematical checks and error discussion; its written form alone does not ensure validity. These are design claims, not measured effects; Hypothesis~\ref{h:chain} states them in testable form.

% ======================================================================
\section{The progression from Grade~1 to Grade~8}\label{sec:progression}
% ======================================================================

The progression is a spiral: every grade has its own chapter, and every chapter revisits earlier content before extending it. Revisiting is implemented literally, by reuse: a later chapter can include an earlier definition, example or exercise as it stands, so that Grade~4, for instance, opens its section on equivalence with the Grade~3 definition and situation. Grades~6 and~7 are consolidation chapters, which cover the whole subject in fourteen sections each before Grade~8 transposes the same structure to algebraic fractions. Each chapter consists of a lesson and an exercise set, and each exercise set ends with written multi-step problems. Table~\ref{tab:progression} summarises the progression; the paragraphs below describe it.

\begin{table}[ht]
\centering\small
\begin{tabularx}{\textwidth}{c>{\raggedright\arraybackslash}p{1.25cm}>{\raggedright\arraybackslash}X>{\raggedright\arraybackslash}p{3.3cm}}
\toprule
Grade & Class by age & New content & Meaning or model introduced \\
\midrule
1 & CP & fair and unfair cuts; halves and quarters of shapes and collections; the notation $\frac{a}{b}$ & fair share of one unit (sharing) \\
2 & CE1 & unit fractions; every fraction as unit fractions; comparison (same denominator, unit fractions); fraction of a collection & partitive reading: $a$ iterations of $\frac{1}{b}$ \\
3 & CE2 & number line and measurement; equivalence (rule by subdivision); comparison on the line; addition and subtraction (same denominator; general procedure through a common multiple); the chain of equalities, without reasons & the bridge: length, point and numerical measure \\
4 & CM1 & fractions greater than one on the line; decimal fractions (tenths, hundredths, decimal point) & counting beyond the unit \\
5 & CM2 & simplification; ordering by common denominator; unlike denominators revisited; models of multiplication and division & area model; ``of'' as $\times$ \\
6 & 6e & consolidation; cross-multiplication; quotient $a\div b=\frac{a}{b}$; operator $\frac{a}{b}\times N$; decimals by division; mixed numbers & quotient (operative definition); operator \\
7 & 5e & rational numbers; multiplication and division of fractions; signs; order of operations & rational number ($b\neq 0$) \\
8 & 4e & algebraic fractions: the Grade~7 structure with letters; domains & domain of an expression \\
\bottomrule
\end{tabularx}
\caption{Overview of the progression. The nominal class is the French class that the level corresponds to by age; Section~\ref{sec:implementation} shows how the French-curriculum chapters actually draw on the levels.}
\label{tab:progression}
\end{table}

\paragraph{Grade 1.} Two parts: halves and quarters (the fairness of a cut; $\frac{1}{2}$ and $\frac{1}{4}$ on discs, rectangles, strips and collections; ``two halves make one whole'', ``a half is larger than a quarter'', in words), then the general notation through one cake cut into equal parts and the definition of Section~\ref{sec:p1}. No operations; denominators up to~$10$ in the exercises.

\paragraph{Grade 2.} The same opening situation and definition; reading fractions aloud; the unit fraction and the decomposition of every fraction into unit fractions; comparison with the same denominator ($\frac{3}{10}<\frac{7}{10}$) and between unit fractions ($\frac{1}{6}<\frac{1}{3}$); fractions of collections by equal sharing ($\frac{1}{3}$ of $12$ beads); first sums and differences with the same denominator in the exercises ($\frac{2}{5}+\frac{3}{5}=1$), which the lesson does not yet define.

\paragraph{Grade 3.} The central level of the elementary part, where the bridge is explicitly taught: the number line (the three-step placement, ``the fraction is a length'', measuring and drawing segments with rulers graduated in fractions of a unit, lengths such as $1+\frac{3}{10}$ beyond~$1$); equivalence by subdivision, with a cake cut into three and then six parts, and $\frac{1}{2}$ and $\frac{2}{4}$ at the same point; comparison on the line; addition and subtraction with the same denominator ($\frac{3}{3}=1$, $\frac{0}{3}=0$) and, for unlike denominators, the procedure through a common multiple, practised when one denominator is a multiple of the other; the chain of equalities under the bars.

\paragraph{Grade 4.} Beyond the unit: an opening situation with $\frac{3}{2}$, fractions greater than one on the line ($\frac{5}{3}=1+\frac{2}{3}$, $\frac{9}{4}$, $3+\frac{1}{2}$), the definition completed (whole numbers; the denominator is never zero); equivalence and addition revisited through the Grade~3 material; decimal fractions (tenths and hundredths, $\frac{47}{10}=4+\frac{7}{10}$, the decimal point, $\frac{7}{10}=0.7$).

\paragraph{Grade 5.} Simplification; the general ordering procedure through a common denominator ($\frac{3}{4}$ and $\frac{5}{8}$); addition and subtraction with any denominators, through the Grade~3 procedure; and the models of multiplication and division of Section~\ref{sec:p2}---a whole number times a fraction on the line, a fraction of a fraction on an area, the number of halves in~$3$, and the prediction of whether a result is larger or smaller.

\paragraph{Grade 6.} Grade~6 consolidates the subject in fourteen sections and adds the quotient (two cakes shared among three people; $a\div b=\frac{a}{b}$, Section~\ref{sec:p4}), the general operator ($\frac{a}{b}\times N=\frac{aN}{b}$), cross-multiplication as a criterion, the conversion of fractions to decimals by division and back by place value---with the question whether a fraction is a whole number, a decimal number that is not whole, or neither ($\frac{12}{4}$, $\frac{3}{4}$, $\frac{4}{3}$), an objective that the French \emph{sixi\`eme} programme also sets \citep{men2025cycle3}---and the relation of fractions greater than one to division with remainder, to the product $n\times\frac{1}{b}$ and to mixed numbers. Ordering uses $1$ and $\frac{1}{2}$ as benchmarks, and the written problems include a historical text on ``broken numbers''.

\paragraph{Grade 7.} Grade~7 places fractions among the rational numbers, the numbers $\frac{a}{b}$ with $a,b$ integers and $b\neq 0$, and completes the arithmetic: simplification by the greatest common factor, multiplication (the area model, the rule, cancelling before multiplying), division through the reciprocal, sign conventions and the comparison of signed fractions, and the order of operations, with the fraction bar as a grouping symbol. Its written problems include problems in the style of the French end-of-collège examination.

\paragraph{Grade 8.} Grade~8, \emph{Algebraic fractions}, keeps the order of the Grade~7 sections (thirteen of them, without fractions greater than one); each section opens with a short text, about half of which recall the numerical case, before one to three worked examples. Pupils simplify expressions such as $\frac{6x}{9}$ and $\frac{9x^3}{30x}$, write algebraic fractions from verbal descriptions, perform the four operations and solve equations with a single fraction ($\frac{12}{4}=\frac{x}{6}$). The domain of every cancellation and the contrast between terms and factors (Sections~\ref{sec:p4} and~\ref{sec:misconceptions}) are new; the section on ratio and operator mentions percentages, discounts and proportional reasoning with letters.

% ======================================================================
\section{Anticipated misconceptions}\label{sec:misconceptions}
% ======================================================================

An a priori analysis states which difficulties a design anticipates and how it responds to them. Table~\ref{tab:misconceptions} lists the misconceptions that the progression addresses, all but the last documented in the literature, together with the design feature intended to counter each and the grade where it appears. The responses are of two kinds: \emph{preventive}, when a representation makes the erroneous reasoning produce a visibly wrong result, and \emph{explicit}, when the chapter confronts the error with a statement or a counter-example.

\begin{table}[ht]
\centering\small
\begin{tabularx}{\textwidth}{>{\raggedright\arraybackslash}p{3.5cm}>{\raggedright\arraybackslash}X>{\raggedright\arraybackslash}p{1.5cm}}
\toprule
Misconception & Design response & Grade \\
\midrule
A fraction is two whole numbers \citep{stafylidou2004,ni2005} & preventive: one bar, one length; one point on the number line; ``the fraction is a length'' & 1, 3 \\
A larger denominator gives a larger fraction \citep{ni2005} & explicit: unit fractions compared (the more parts, the smaller each part); same numerator on the line & 2, 3, 6 \\
Adding numerators and denominators separately \citep{streefland1991,ni2005} & preventive: only parts of the same size are counted together; re-partitioning before adding & 3, 5 \\
Multiplication makes bigger, division makes smaller \citep{fischbein1985} & explicit: predicting whether a result is larger or smaller; models of $3\times\frac{2}{3}$ and $3\div\frac{1}{2}$ & 5 \\
Equivalent fractions are different numbers & preventive and explicit: the same point on the line; two ways of writing the same number & 3, 6 \\
A part--part ratio read as a part of the whole & explicit, in the ratio strand: what the fraction measures in each kind of ratio & 4 (ratio) \\
The fraction known as an action (``$\frac{2}{3}$ of'') rather than as a number & preventive: the number line before the operator; ``of'' as a sharing situation in the ratio strand & 3, 6 \\
Cancelling terms of a sum \citep{kirshner2004} & explicit: $\frac{x+2}{2}$ and $x$ are not identically equal ($x=0$ gives $1\neq 0$), against $\frac{2x+4}{2}=x+2$ & 8 \\
Losing the domain after a simplification & explicit: the excluded value written next to every cancellation in the lessons & 8 \\
\bottomrule
\end{tabularx}
\caption{Anticipated misconceptions and the design responses of the progression. The last row is a design observation rather than a documented finding.}
\label{tab:misconceptions}
\end{table}

\paragraph{The fraction as one number.} Pupils who see a fraction as two whole numbers separated by a bar compare and combine the two numbers separately \citep{stafylidou2004}. The response is preventive and runs through the whole progression: the definition is read on one bar, which is one length, and from Grade~3 the fraction is one point on the number line. Placing $\frac{2}{3}$ between $0$ and $1$ is intended to establish its relation to other numbers; successful placement alone does not establish the pupil's conception.

\paragraph{Whole-number bias.} Three errors apply whole-number properties to fractions \citep{ni2005}; \citet{streefland1991} calls such errors ``N-distractors''. That a larger denominator gives a larger fraction is confronted explicitly from Grade~2, where unit fractions are compared, and on the number line in Grade~3. Adding numerators and denominators separately is meant to be prevented by the linear model (Figure~\ref{fig:operations}): only parts of the same size can be counted together, and the Grade~3 situation on $\frac{1}{2}+\frac{1}{4}$ states that the numerators cannot be added directly. The progression does not, however, confront the erroneous sum $\frac{1}{2}+\frac{1}{3}=\frac{2}{5}$ with a counter-example, although errors of this kind are frequent and persistent \citep{lortieforgues2015}. A design that wanted an explicit response would add one---for instance by observing that $\frac{2}{5}<\frac{1}{2}$, so that the ``sum'' would be smaller than one of its terms---and would pair it with the ratio situation of Section~\ref{sec:p3} in which adding numerators and denominators is correct. That multiplication always makes bigger and division smaller \citep{fischbein1985} is addressed in Grade~5, where pupils predict the size of a result before computing it.

\paragraph{Referent and equivalence.} That equivalent fractions name different numbers is countered by the same-point representation in Grade~3 and by the statement, in Grade~6, that equal fractions are two ways of writing the same number. The reading of a part--part ratio as a part of the whole is the object of the explicit distinction of the ratio strand (Section~\ref{sec:p3}).

\paragraph{Algebraic fractions.} Two errors appear once letters enter. Pupils cancel a term that appears in the numerator and in the denominator, reading $\frac{x+2}{2}$ as $x$, an error that \citet{kirshner2004} relate to the visual salience of algebraic notation and to the unstated distinction between terms and factors; the Grade~8 lesson contrasts this with a true common factor. And pupils who simplify correctly lose the excluded value: $\frac{9x^3}{30x}$ and $\frac{3x^2}{10}$ coincide everywhere except at $x=0$, where the first is undefined. Writing the excluded value next to every cancellation is the response. In the current exercise set, however, only three of the thirteen simplification statements with a letter in the denominator give the condition $x\neq 0$, and no answer restates it, although the lesson does; we return to this in Section~\ref{sec:discussion}.

\subsection{A task-level example: placing a fraction on the number line}\label{sec:taskvariables}

The analysis above stays at the level of the design. An a priori analysis in the full sense also examines tasks: which \emph{didactic variables} a task offers (Artigue speaks of \emph{variables de commande}), that is, which of its features the designer can choose, and which strategies and errors each choice makes likely \citep{artigue1990}. Placing $\frac{a}{b}$ on a number line, the central task of Principle~\ref{pr:linear}, has at least four such variables.
\begin{enumerate}[label=(V\arabic*)]
  \item \emph{The partition given.} The unit segment may be already cut into $b$ equal parts, cut into a multiple of $b$, or not cut at all, in which case the pupil partitions it. A given partition into $b$ parts reduces the task to counting; it also allows a pupil to count tick marks rather than intervals, and so to place $\frac{3}{8}$ one graduation too early or too late.
  \item \emph{The extent of the line.} If the drawn line runs from $0$ to $1$, the unit and the drawn line coincide; if it runs further, a pupil who reads $\frac{a}{b}$ as ``$a$ parts out of all the parts drawn'' takes the whole line as the unit, an error documented on number-line tasks \citep{bright1988}.
  \item \emph{The relation between graduation and denominator.} Placing $\frac{1}{2}$ on a line graduated in quarters requires recognising $\frac{1}{2}=\frac{2}{4}$, which turns a placement task into an equivalence task.
  \item \emph{The size of the fraction.} A target greater than one cannot be reached by counting parts ``out of'' the unit, and separates the counting reading from the iterative one.
\end{enumerate}
The exercises of the progression move along these variables. In Grade~3, pupils first place fractions on a line from $0$ to $1$ whose unit is already cut into the right number of parts (V1 given, V2 and V4 at their simplest), then graduate the unit themselves (V1 free); the equivalence exercises ask for a fraction on a line graduated in a multiple of its denominator, such as $\frac{2}{4}$ on a unit cut into eighths (V3), and the measurement exercises use rulers graduated in quarters, fifths or tenths that extend beyond one unit, with lengths such as $2+\frac{1}{4}$ (V2). Grade~4 then asks for fractions greater than one (V4). Hypotheses~\ref{h:improper} and~\ref{h:magnitude} (Section~\ref{sec:discussion}) can be made observable through these variables: for instance, the proportion of pupils who place $\frac{5}{3}$ correctly on a line from $0$ to $2$ graduated in thirds (V2 and V4 at their most demanding) is an indicator---not a proof, since a counting routine or a conversion to a mixed number can also succeed---of the iterative reading.

% ======================================================================
\section{Implementation}\label{sec:implementation}
% ======================================================================

The progression exists as sixteen bilingual (English and French) source documents, a lesson and an exercise set for each grade, written in \LaTeX\ and stored in a content repository together with the companion ratio strand. A curriculum-specific chapter is assembled from these sources by reference: it lists, in its own classroom order, blocks of the source progression (definitions, examples, exercises, videos), so that one block can appear in several curricula without being copied, and a correction made in the source reaches every chapter that uses it. The curriculum chapters are published as web pages with interactive exercises, and as printable documents, at \url{https://commeunjeu.org}. The source documents themselves are not published: what is public is the curriculum chapters assembled from them.

Two curricula currently draw on the fraction progression: the French national curriculum, from CE1 to \emph{troisi\`eme}, and the International Baccalaureate, from the first year of the Primary Years Programme to the fourth year of the Middle Years Programme. Appendix~\ref{app:resources} lists the public chapter addresses. The French chapters do not follow the levels one to one (Table~\ref{tab:frmap}): the Grade~1 material is taught in CE1, because the French programmes in force since 2025 introduce fractions in CE1 \citep{men2024cycle2};
Grade~6 is the main source of the CM2, \emph{sixi\`eme} and \emph{cinqui\`eme} chapters, and Grade~7 of the \emph{quatri\`eme} and \emph{troisi\`eme} chapters. Grade~8, on algebraic fractions, is the main source of two IB chapters, which use 127 of its 166 blocks (answers excepted); in the French curriculum it supplies only part of the chapter on literal calculus of \emph{troisi\`eme}. Conversely, some sections of the progression are used by no public fraction chapter: the cross-multiplication sections of Grades~6 and~7, the main example of Principle~\ref{pr:technique}, are among them. The analysis can therefore be checked against the public chapters for most, but not all, of what it describes. This decoupling of a source progression from its curriculum projections is what allows a single didactic design to serve curricula whose year-by-year allocation differs; it also means that a projection can leave out part of the design.

\begin{table}[ht]
\centering\small
\begin{tabular}{llll}
\toprule
French class & Chapter (address in Appendix~\ref{app:resources}) & Main source & Other sources \\
\midrule
CE1 & First fractions & Grades 2, 1 & 3, 4, 5 \\
CE2 & Fractions & Grade 3 & 1, 2, 4, 5 \\
CM1 & Fractions & Grade 4 & 1, 2, 3, 5, 6 \\
CM2 & Fractions & Grade 6 & 3, 4, 5, ratio strand \\
6e & Fractions as numbers & Grade 6 & 3, 4, 5 \\
6e & Calculating with fractions & Grade 6 & 2, 3, 4, 5 \\
5e & Rational numbers & Grade 6 & 2, 3, 4, 5 \\
4e & Rational numbers & Grade 7 & 3, 4, 5, 6 \\
3e & Rational numbers & Grade 7 & 3, 4, 5, 6 \\
\bottomrule
\end{tabular}
\caption{Levels of the progression from which each French-curriculum chapter draws its blocks (main source: the level that supplies the largest number of blocks).}
\label{tab:frmap}
\end{table}

\paragraph{Source order and experienced order.} The source levels describe the designer's organisation, while curriculum assembly determines which blocks are presented together and in what order. Neither establishes the sequence actually taught. A claim that numerical meaning precedes explicit operator treatment for a given group would require tracing the assembled fraction and ratio chapters together, including prerequisite references, omitted blocks and teacher choices. Table~\ref{tab:frmap} records provenance rather than that chronology. The present analysis therefore identifies intended sequencing and documentary tensions; it does not report observed pupil trajectories.

Table~\ref{tab:trace} records, for each principle, the blocks of the corpus that carry it; it is the table against which Section~\ref{sec:discussion} was written.

\begin{table}[ht]
\centering\small
\begin{tabularx}{\textwidth}{>{\raggedright\arraybackslash}p{2.6cm}>{\raggedright\arraybackslash}X>{\raggedright\arraybackslash}X}
\toprule
Principle & Lessons (grade: block) & Exercises (grade: section) \\
\midrule
1 Share of one unit & 1: opening problem (fair cut), definition ``Fraction''; 2: definition ``Unit fraction'', proposition ``Every fraction is made of unit fractions''; 4: definition completed (whole numbers, denominator never zero); 6: proposition $a\div b=\frac{a}{b}$ & 1: equal parts, fractions in diagrams, drawing; 2: building from unit fractions \\
2 The bridge & 3: methods ``Placing a fraction'', ``Measuring with a ruler'', ``Drawing a segment'', text ``The fraction is a length''; 4: method ``Placing a fraction greater than~1'' & 3: graduating the unit strip, reading a length, bounding, drawing; 4: number line greater than~1 \\
3 Numbers before operators & Ratio 3: method ``Finding a part in part--whole relationships''; Ratio 4: definitions part--part and part--whole; Ratio 5: method ``From `half of' to `multiply by'\,''; Fraction 6: opening problem and proposition $\frac{a}{b}\times N$ & Ratio 3--5: part in part--whole; Fraction 6: ratio and operator \\
4 Technique before formalism & 3: definition ``Equivalent fractions''; 6: proposition and remark on cross-multiplication; 7: proposition with $b,d\neq 0$; 8: examples with ``valid for $x\neq 0$'' & 6: cross-multiplication; 8: simplification (condition in 3 of the 13 statements with a letter in the denominator) \\
5 Chain of equalities & 3: example $\frac{1}{2}+\frac{1}{4}$ (tabular); 4: example $\frac{9}{4}$ with reasons. Chains with a reason: 4: 6 of 6; 5: 5 of 5; 6: 6 of 9; 7: 14 of 22; 8: 11 of 12 & Chains with a reason: 2: 0 of 8 (two-line chains); 4: 4 of 48; 5: 9 of 66; 6: 15 of 106; 7: 44 of 164; 8: 69 of 93 \\
\bottomrule
\end{tabularx}
\caption{Traceability of the principles to the blocks of the corpus (lesson blocks named by their titles; counts exhaustive over the English half of the mark-up, with the criteria of Section~\ref{sec:intro}).}
\label{tab:trace}
\end{table}

Twenty-eight videos accompany the lessons, one per section for most sections of Grades~1 to~7 and one for the first section of Grade~8; twenty-two exist in English and in French, and six in French only. Each is attached to one lesson section and embedded in the chapters that use that section (Appendix~\ref{app:resources}).

% ======================================================================
\section{Discussion and limitations}\label{sec:discussion}
% ======================================================================

\subsection{Position with respect to existing recommendations}

Several features of the progression coincide with existing recommendations. The reading of $\frac{a}{b}$ as $a$ parts of size $\frac{1}{b}$, the placement of fractions on the number line, and equivalence understood as ``the same point'' follow the Common Core standards for Grades~3 and~4 \citep{ccss2010}, and the central role of the number line agrees with the practice guide of \citet{siegler2010}. The progression departs from these references in three ways, each of which is a hypothesis rather than an established advantage.

\emph{An earlier start and a denser third year.} The notation $\frac{a}{b}$ appears at the first level of the progression (nominally age~6), with a restricted scope (halves, quarters and simple fractions, no operations), one year before the French programmes, which start in CE1 \citep{men2024cycle2}, and two years before the Common Core standards, which limit Grades~1 and~2 to partitioning shapes into halves, thirds and fourths and introduce the notation in Grade~3.
Grade~3 gathers the number line, equivalence, comparison and the first additions, which the Common Core standards distribute across Grades~3--5. The intended benefit of the earlier levels is familiarity with fractional notation before number-line work; the risk is that the notation remains a label for pictures. The shared representation offers an organising principle: equivalence can be read as equal length, comparison as longer or shorter, and addition as concatenation. These are nevertheless distinct learning problems. In Grade~3, pupils are expected to place and measure fractions, recognise simple equivalences, and add by repartitioning when one denominator is a multiple of the other. The general common-multiple procedure is stated, while practice with arbitrary denominators belongs to Grade~5. Understanding the reference unit, conserving a measure under subdivision and coordinating an operation with its diagram cannot be reduced to a single insight that ``the fraction is a length''. Whether the common representation supports manageable pacing is an empirical question; the design proposes a sequence, without demonstrating its cognitive cost.

\emph{A later quotient.} The quotient meaning comes one year later than in the Common Core standards, where it belongs to Grade~5, and the general multiplication rule two years later. The models of multiplication appear in the progression's Grade~5, one year after the standards' product of a fraction by a whole number (Grade~4) and in the same year as their area model; decimal fractions, on the other hand, appear in Grade~4 in both. The progression thus lets the additive structure and the models stabilise before the multiplicative meanings are stated in general. With respect to the French programmes, the progression agrees on the quotient in \emph{sixi\`eme}, on the number line as what gives fractions the status of numbers, and on fractions greater than one in CM1 \citep{men2025cycle3}; it differs in starting one year earlier; in stating, in its Grade~3, the general procedure for unlike denominators through a common multiple, where the CE2 programme asks only for the case in which one denominator is a multiple of the other and keeps fractions at most one with denominators up to twelve \citep{men2024cycle2}; and in the status it gives the operator (Section~\ref{sec:p3}).

\emph{Sharing as a later situation.} Equal sharing of several objects among several people, which the practice guide recommends building on early \citep{siegler2010} and which research programmes in the tradition of realistic mathematics education and cognitively guided instruction use as the entry into fractions \citep{streefland1991,empson2011}, gives meaning to the quotient only in Grade~6. Before that, sharing is the partition of a single object, or the equal sharing of a collection to find a unit fraction of it. Whether this delay costs something that the early number line does not compensate is an open empirical question.

The contrast with the engineering of \citet{brousseau2004}, and with the commensuration method studied by \citet{ratsimba1981}, is of a different kind. There, the measure of an object is not obtained by partitioning a given unit: the need for a new number arises from a measurement problem that whole numbers cannot solve, and the rational number first appears as a relation between multiples of the object and multiples of a reference. The present progression gives the unit and asks only that it be partitioned. It thereby gives up the epistemological depth of commensuration, in which pupils experience why rational numbers are needed, in favour of an entry intended to limit initial coordinations and used from the first source level. The comparative learning cost has not been measured. A progression could combine the two, for instance by introducing a commensuration situation once the number line is in place; this one does not.

\subsection{What the analysis reveals about the chapters}

Confronting the principles with the chapters as they stand brings out points on which the chapters do not yet carry the design. They are of three kinds.

\emph{Conceptual: statements or justifications whose scope needs clarification.}
\begin{enumerate}[label=(\alph*)]
  \item\label{it:titles} The fraction strand's sections of Grades~6--8 are titled ``Fraction as a ratio and operator'', although the Grade~6 and~7 texts do not discuss ratios. The Grade~8 phrase ``as a ratio, the fraction $\frac{1}{4}$ describes one part out of four equal parts'' also needs to distinguish a part--whole relation from a part--part comparison. Treating operator and ratio as later meanings is compatible with Principle~\ref{pr:ratio}; their mere mention is not a contradiction. The recommended revision is to use a title matching the content, such as ``Taking a fraction of a quantity'', and connect these uses explicitly to numerical meaning.
  \item\label{it:division} Grade~7 introduces division by a fraction through multiplication by the reciprocal and the inverse relation between division and multiplication, without recalling the Grade~5 measurement model ($3\div\frac{1}{2}=6$). An inverse-multiplication argument can be valid; the absence of the earlier picture does not invalidate it. The recommended revision is to make the argument explicit, including the nonzero divisor, and connect it to the earlier measurement interpretation.
  \item\label{it:withdraw} Conditions are introduced but not always withdrawn: Grade~7 reuses a subtraction procedure stating that the first fraction must not be smaller than the second, and then computes $\frac{2}{3}-\frac{3}{4}=-\frac{1}{12}$.
\end{enumerate}
\emph{Sequencing: references that precede their explicit treatment in the sources.}
\begin{enumerate}[label=(\alph*),resume]
  \item\label{it:ratiostrand} The ratio strand's Grade~2 definition writes $a:b$ ``or as the fraction $\frac{a}{b}$'', simplifies ratios ($2:4=1:2$) and uses a part--whole $\frac{2}{4}$. These statements precede the source fraction strand's equivalence lesson and the ratio strand's distinction between part--part and part--whole relations. The Grade~4 quotient remark provides a further early encounter that the Grade~6 quotient section does not revisit. Coordinating these passages requires decisions about content order across the two strands, not only changing labels. Their effect on a curriculum projection must be checked in its actual assembled order.
  \item\label{it:forward} Two forward references break the order of a chapter: the Grade~3 method for recognising equivalent fractions uses ``the same point on the line'' one section before the number line is introduced, and Grade~6 says ``use the quotient'' before its section on the quotient.
  \item\label{it:crossmult} The Grade~6 cross-multiplication remark uses cancellation before the later quotient section, without developing the available justification by iteration or subdivision (Section~\ref{sec:p4}). The issue is an unstated local justification, not unavoidable mathematical circularity.
\end{enumerate}
\emph{Editorial: wording or practice that falls short of the principle.}
\begin{enumerate}[label=(\alph*),resume]
  \item\label{it:whole} The definition takes \emph{the unit} as its referent, but the Grade~1 and Grade~2 texts on halves and unit fractions speak of ``the whole'', and so does the Grade~2 answer quoted in Section~\ref{sec:p5}.
  \item\label{it:answers} The Grade~8 lessons write the excluded value next to every cancellation; in the exercises, three of the thirteen simplification statements with a letter in the denominator give the condition, and no answer restates it. The same holds of Principle~\ref{pr:chain} in the middle grades: the lessons annotate most of their chains from Grade~4, whereas the exercise answers annotate fewer than one chain in six before Grade~7, about one in four in Grade~7 and three in four in Grade~8 (Table~\ref{tab:trace}). The lessons and answers thus model the intended habits to different degrees. Omitting a repeated domain condition is not automatically an invalid answer when it is inherited from the question. The instructions must also be examined before inferring what pupils are required to justify.
  \item\label{it:counter} No chapter confronts the sum $\frac{1}{2}+\frac{1}{3}=\frac{2}{5}$ with an explicit counter-example, although it is a classic whole-number error on addition.
\end{enumerate}
These findings have different consequences. Terminology and missing annotations can often be revised within existing blocks. Forward references and the relation between ratio and fraction definitions may require reordering or rewriting content across strands. Material validation of measurements, the subdivision proof and a general quotient argument are proposed additions, not documented features of the current chapters. The principles can be retained as intentions while their realisation changes substantially; it would be misleading to describe every discrepancy as editorial. Task~B is different again: its absence from the initial sequence is deliberate and remains part of the design to be evaluated. Table~\ref{tab:status} separates observed content, limits and recommended revisions.

\begin{table}[ht]
\centering\footnotesize
\begin{tabularx}{\textwidth}{>{\raggedright\arraybackslash}p{1.9cm}>{\raggedright\arraybackslash}p{1.5cm}>{\raggedright\arraybackslash}X>{\raggedright\arraybackslash}X>{\raggedright\arraybackslash}X}
\toprule
Principle & Status & Realised in the chapters & Observed limits & Recommended change \\
\midrule
1 Share of one unit & present; extensions to justify & definition on a bar (G1); unit fractions (G2); beyond one (G4); quotient (G6) & ``whole'' for ``unit'' in G1--2 (item~\ref{it:whole}); quotient stated from one example (Section~\ref{sec:p4}) & write ``unit'' throughout; justify $a\div b=\frac{a}{b}$ in general in G6 \\
2 The bridge & tasks present; learning untested & ruler, drawing and number line (G3); beyond one (G4); stable reference & ``same point on the line'' before the line (item~\ref{it:forward}); external validation of readings (Section~\ref{sec:crossing}); no counter-example to $\frac{1}{2}+\frac{1}{3}=\frac{2}{5}$ (item~\ref{it:counter}) & reorder the G3 chapter; add a check of each reading (drawing it; counting pieces); add the counter-example \\
3 Numbers before operators & partly applied & operator written in the fraction strand only in G6; ``of'' as sharing in the ratio strand (G3, G5) & titles do not match content; early fraction notation in G2 ratios; source order is not pupil chronology (items~\ref{it:titles}, \ref{it:ratiostrand}) & clarify uses and referents; coordinate both strands and defer the general ratio--fraction link; revisit the G4 remark in G6 \\
4 Technique before formalism & partly applied & equivalence by subdivision (G3); conditions stated (G4, G7; G8 lessons) & local arguments need development (items~\ref{it:crossmult}, \ref{it:division}); subtraction condition retained after signed numbers (item~\ref{it:withdraw}) & add subdivision proof; explain reciprocal division and its domain; update subtraction for signed numbers \\
5 Chain of equalities & lessons: mostly applied; answers: partly & chains from G2 (answers) and G3 (lesson); reasons in most lesson chains from G4 & few annotated answers before G8; excluded value absent from G8 answers (item~\ref{it:answers}) & annotate model answers from G4; restate the condition in G8 answers \\
\bottomrule
\end{tabularx}
\caption{Status of each principle in the chapters (version of 9~October 2026): what they do, what the analysis observes, and what it recommends. Statuses concern the dated source documents, not demonstrated learning. The recommended changes are proposals arising from this audit, not changes reported as implemented in that corpus.}
\label{tab:status}
\end{table} These points illustrate the use of an explicit design analysis: it gives a reader, and the author, a standard against which the materials can be checked.

\subsection{What the progression does not address}

Two aspects of rational numbers remain outside the main scope of the eight source levels. The first is \emph{density}: between two distinct rational numbers there is another. Pupils may believe that no fraction lies between $\frac{2}{5}$ and $\frac{3}{5}$ \citep{vamvakoussi2004}. The broader design places the formal study of density in Grade~10, with number sets and intervals, beyond the corpus analysed here. That later formalisation should be distinguished from earlier examples of inserting a fraction between two others: $\frac{2}{5}=\frac{4}{10}<\frac{5}{10}<\frac{6}{10}=\frac{3}{5}$. Subdivision supplies the mathematical resources for such examples, but the eight fraction chapters do not make this an explicit objective. The companion decimal strand includes insertion between two decimals in Grade~6. This does not establish transfer to fractions, and the possible persistence of a discrete view is a limitation to investigate. The second aspect is the relation between fractions and decimal expansions beyond conversion techniques. Decimal fractions appear in Grade~4 and terminating and non-terminating decimals in Grade~6; a repeated remainder is observed in one long division. Periodicity and the reverse question of which decimal expansions represent fractions are left to the decimal strand and are not evaluated here.

\subsection{Hypotheses for an empirical evaluation}

Because no learning data are reported, the claims of Sections~\ref{sec:p1}--\ref{sec:p5} are formulated as hypotheses. ``Grade~$N$'' refers to the source level used, not necessarily pupils' school class. A study would need to document the actual sequence across both strands and the activities assigned, rather than treating source labels as exposure measures. Primary hypotheses concern the bridge; secondary ones concern supporting decisions. Comparisons should match time, content coverage and assessment conditions. To test the bridge specifically, both conditions should include number-line work, with the explicit coordination of sharing, measurement and points being the intended contrast. Predictions of ``no loss'' require a prespecified educationally meaningful non-inferiority margin and adequate precision; a nonsignificant difference alone would not establish equivalence.

\paragraph{Primary.}

\begin{hypothesis}\label{h:magnitude}
Explicit coordination of a shared bar, its measured length and a point on the number line improves transfer between representations by Grade~5, beyond practice with the representations separately. Measures include unpartitioned number-line estimation on $[0,1]$ and $[0,2]$, and comparisons with different numerators and denominators. Transfer items ask pupils to recognise or produce the same fraction across a partitioned bar, a length relative to a supplied unit and an unpartitioned line, without prescribing a conversion procedure. Diagram size, orientation and the extent of the drawn line vary while the reference is specified. A combined share--measure--place item assesses coordination with support; separate, unfamiliar items assess less-cued transfer. Code identification of the unit, agreement of numerical values, explanations and whole-number-bias errors separately, so consistent but incorrect responses do not count as successful transfer.
\end{hypothesis}

\begin{hypothesis}\label{h:units}
Keeping the reference stable in the initial fraction tasks supports numerical meaning without a practically important later disadvantage in coordinating units. In Grades~3--5, measures include locating fractions on lines extending beyond~$1$, measuring bars longer than the designated strip and recognising the same fraction in proportionally resized diagrams. At the end of Grade~6, use two distinct probes: (i) two strips, $20$~cm and $10$~cm, each designated as one unit, whose halves have equal numerical measures in their respective units but different physical lengths; (ii) one unchanged $6$~cm segment measured with reference strips of $8$~cm and $4$~cm, giving $\frac{3}{4}$ and $\frac{3}{2}$. Record whether pupils preserve the specified reference and distinguish number from physical length. The second probe tests reunitisation directly; operator success is not its substitute. Compare with a curriculum introducing unit changes earlier, record any explicit teaching of reunitisation in either group, and distinguish unprompted transfer from performance after such instruction. Non-inferiority is a testable prediction, not a conclusion licensed by the present design analysis.
\end{hypothesis}

\begin{hypothesis}\label{h:operator}
Prior work establishing fraction magnitude, followed by explicit treatment of scalar multiplication, improves pupils' ability to use $\frac{2}{3}$ as a number without an immediate part--whole context, while preserving performance on applications to quantities. At the end of Grade~6, assess context-free placement and comparison alongside ``$\frac{2}{3}$ of $30$'', with matched comparison groups. Ask what $\frac{2}{3}$ means and allow words or drawings. Code action descriptions, perceptual part--whole descriptions, measures relative to a unit and numerical relations as distinct response categories; they are not assumed to form a validated developmental scale. Interpret them with performance and explanations across several items, rather than treating one preferred wording as evidence of understanding. Record the actual timing of explicit operator instruction in both strands.
\end{hypothesis}

\paragraph{Secondary.}

\begin{hypothesis}\label{h:improper}
Pupils who meet unit-fraction iteration in the Grade~2 material (``three eighths is three times one eighth'') place fractions greater than one more accurately when they reach the Grade~4 material than pupils taught the counting reading alone. Measure: placing $\frac{5}{3}$ and $\frac{7}{4}$ on lines from $0$ to $2$ that are either unpartitioned or partitioned in a non-matching way (sixths, eighths), so that counting the target numerator's number of given intervals does not suffice; code refusals, placements below~$1$ and use of the whole line as the unit \citep{bright1988}. Ask pupils to identify the unit fraction they iterated. Correct placement remains an indicator rather than proof of an iterative scheme, since learned conversions can also succeed.
\end{hypothesis}

\begin{hypothesis}\label{h:ratio}
The explicit part--part/part--whole distinction of the ratio strand (Principle~\ref{pr:ratio}) reduces the reading of a part--part ratio $a:b$ as ``$\frac{a}{b}$ of the total''. Measure: problems giving a part--part ratio and asking for the fraction of the whole, at the end of Grade~6.
\end{hypothesis}

\begin{hypothesis}\label{h:chain}
Documented practice in explaining successive fraction transformations from the primary levels reduces invalid transformations and improves error detection in later algebraic simplifications. Compare with instruction matched for mathematical content and opportunities to explain reasoning, while varying the use of written equality chains. Assess unfamiliar expressions, valid and invalid transformations, and reasons for accepting or rejecting them. Record use of the taught layout separately as compliance with a convention. The current corpus's uneven annotations mean that this hypothesis concerns the intended practice; exposure must be checked in an empirical study.
\end{hypothesis}

\begin{hypothesis}\label{h:domain}
Explicitly retaining the original domain beside algebraic cancellations reduces domain loss after simplification in Grade~8. Items ask for the original domain, a simplified expression on that domain and whether both expressions have the same natural domain. Correctly inheriting a condition from the question is distinguished from asserting an unrestricted identity.
\end{hypothesis}

\begin{hypothesis}\label{h:cancel}
The explicit contrast between $\frac{x+2}{2}$ (not identically $x$) and $\frac{2x+4}{2}$ in Grade~8 reduces the cancellation of terms of a sum. Measure: simplification items mixing sums and products in numerator and denominator.
\end{hypothesis}

Answers collected by the online exercises (Section~\ref{sec:implementation}) could give a first observational view of some of these questions, but such data have no comparison group and come from a self-selected population; they could suggest where to look, not establish an effect.

\subsection{Limitations}

The progression was designed by a single author, and the analysis in this paper is that author's. The literature on which it relies is mainly English-language cognitive and educational psychology and the French didactic tradition; other traditions, such as the Japanese and Chinese curricula, in which the teaching of fractions has been studied extensively, are not discussed. Finally, a progression is not a lesson: how teachers enact the chapters, which exercises they choose and how they discuss errors are not described by the documents analysed here. No data on teachers' actual practice with these chapters are available. Which exercises are assigned, how the bar and the line are discussed in class, and whether the reasons that the lessons write are required of pupils are unknown; the online platform records pupils' answers, not teaching. The analysis is also bound to one version of a corpus that changes: several points of Table~\ref{tab:status} may be corrected in later versions of the chapters. The unpublished corpus and counting script limit independent reproduction of the counts and inspection of block-level judgements. A dated archive with stable block identifiers and the coding record would strengthen that evidence; the present revision of the article does not supply such an archive or constitute a new audit of the teaching sources.

% ======================================================================
\section{Conclusion}\label{sec:conclusion}
% ======================================================================

This paper analyses a proposed progression from sharing one unit to treating fractions as numbers with magnitude, and then to explicit operator and quotient interpretations. The bar offers mathematical continuity between a shared object and a measured length; pupils' coordination of those readings remains a learning problem. Keeping the initial reference stable is intended to focus attention on order, equivalence and arithmetic. It postpones reunitisation rather than denying its legitimacy, and neither its cognitive benefit nor its longer-term cost is established here. The subdivision proof of cross-multiplication illustrates how a valid local justification can precede formal generalisation. Equality chains provide a possible written support whose value must be distinguished from adherence to its format. The audit identifies realised features, discrepancies requiring local or curricular revision, and additions proposed by the analysis. Source levels and assembled chapters are not evidence of the order experienced by pupils. The contribution is therefore an explicit, criticisable organisation of meanings and design decisions, together with testable hypotheses and a documented account of the current limits of its realisation. Empirical study and a publicly inspectable source archive remain necessary to evaluate learning effects and independently reproduce the corpus analysis.

\section*{Competing interests}
The author is the founder of CommeUnJeu, which publishes the materials analysed in this paper and offers paid services on its website.

\section*{Acknowledgements}
The author used AI language tools while preparing this paper: Claude (Anthropic) assisted with drafting, with the conversion of an earlier teacher's guide into the present article and with the checking of bibliographic references, and ChatGPT (OpenAI) assisted with critical reading and subsequent revision of the manuscript. The author designed the progression, checked every reference against its source, and takes full responsibility for the content.

\bibliographystyle{abbrvnat}
\bibliography{fractions_refs}

\appendix
% ======================================================================
\section{Public resources}\label{app:resources}
% ======================================================================

The chapters analysed are served at the addresses below (checked on 9 October 2026), in English and in French; the address slugs are in English for both curricula.

\paragraph{Chapters.} All addresses begin with \texttt{https://commeunjeu.org/cours/}; in the French slugs, \texttt{grade-}$N$ denotes the $N$-th year of French compulsory schooling after kindergarten (CP is 1), not a level of the progression.
\begin{center}\small
\begin{tabular}[t]{ll}
\toprule
\multicolumn{2}{l}{French national curriculum} \\
\midrule
CE1 & \texttt{fr-grade-2-first-fractions/} \\
CE2 & \texttt{fr-grade-3-fractions/} \\
CM1 & \texttt{fr-grade-4-fractions/} \\
CM2 & \texttt{fr-grade-5-fractions/} \\
6e & \texttt{fr-grade-6-fractions-as-numbers/} \\
6e & \texttt{fr-grade-6-calculating-with-fractions/} \\
5e & \texttt{fr-grade-7-rational-numbers/} \\
4e & \texttt{fr-grade-8-rational-numbers/} \\
3e & \texttt{fr-grade-9-rational-numbers/} \\
\bottomrule
\end{tabular}\hfill
\begin{tabular}[t]{ll}
\toprule
\multicolumn{2}{l}{International Baccalaureate (PYP, MYP)} \\
\midrule
PYP 1 & \texttt{ib-pyp-1-halves-and-quarters/} \\
PYP 2 & \texttt{ib-pyp-2-fractions/} \\
PYP 3 & \texttt{ib-pyp-3-fractions/} \\
PYP 4 & \texttt{ib-pyp-4-fractions/} \\
PYP 5 & \texttt{ib-pyp-5-fractions-equivalence-and-comparison/} \\
PYP 5 & \texttt{ib-pyp-5-operations-with-fractions/} \\
MYP 1 & \texttt{ib-myp-1-fractions/} \\
MYP 2 & \texttt{ib-myp-2-fractions/} \\
MYP 3 & \texttt{ib-myp-3-rational-numbers/} \\
MYP 3 & \texttt{ib-myp-3-algebraic-expressions-and-fractions/} \\
MYP 4 & \texttt{ib-myp-4-algebraic-fractions/} \\
\bottomrule
\end{tabular}
\end{center}

\paragraph{Videos.} Each chapter page embeds the videos attached to the sections it uses. They are published on two channels, one per language: \texttt{youtube.com/@CommeUnJeuEN} (English) and \texttt{youtube.com/@CommeUnJeu} (French).

\clearpage
\section{Audit conventions and evidence status}\label{app:audit}

This appendix makes explicit the scope of the counts reported in Section~\ref{sec:intro} and Table~\ref{tab:trace}. It does not add observations or certify a new execution of the author's counting script.

\begin{enumerate}[leftmargin=*]
  \item \textbf{Corpus and language.} The reported audit concerns the sixteen fraction source documents and five ratio lessons dated 9~October~2026. Annotation counts use the English source text only; bilingual copies and repeated appearances in assembled curriculum chapters are not additional pupil exposures.
  \item \textbf{Documentary unit.} Presence and source-order claims concern typed blocks. The counts of written chains use multi-line calculations in \texttt{align} or \texttt{aligned} environments. They do not exhaust every mathematically possible equality chain: single-line expressions and tabular picture--calculation displays are outside this particular count. Lessons and model exercise answers are reported separately.
  \item \textbf{Annotation criterion.} A chain is counted as carrying a reason when at least one calculation line has a parenthesised reason, such as ``divide both by $2$''. Words belonging to an expression, such as ``number of girls'', do not count. A reason appearing only in neighbouring prose is excluded from this measure, even when it explains the calculation. The stricter count described in the method requires an annotation on every line except the last result line. These are counts of annotation practice, not judgements that every transformation is valid or justified.
  \item \textbf{Conditions and context.} The appearance of a nonzero-denominator condition in an exercise statement is distinguished from its repetition in the answer. The instructions and surrounding text must be retained when judging whether an answer inherits a restriction. An absent repeated condition is not, by itself, an unrestricted mathematical assertion.
  \item \textbf{Order and instructional status.} Definitions, propositions, methods and remarks provide evidence about the written treatment of a meaning. Their labels do not establish classroom institutionalisation or pupils' prior knowledge. A source-order discrepancy should be checked against the assembled chapter and the other strand before a claim about the order actually taught is made.
\end{enumerate}

For example, a multi-line calculation with a parenthesised reason on one line meets the broad annotation criterion even if its other transformations are unannotated. The same explanation in the paragraph before the calculation does not meet that criterion. This difference is deliberate: the reported statistic describes visible line annotations rather than all available mathematical explanation.

An independently inspectable reanalysis would retain, for each coded item, the dated source file and stable block identifier, language, block type, complete relevant excerpt, chain and annotation decisions, domain context, and reasons for ambiguous classifications. It would distinguish observed content, the analyst's interpretation and a proposed revision. Those item-level records and a fixed source archive are not provided with this manuscript; public curriculum chapters allow partial inspection but cannot substitute for omitted source sections or reproduce all source counts.

\end{document}